\documentclass[a4paper]{amsart}

\usepackage{graphicx}
\usepackage{amsmath}
\usepackage{amssymb}
\usepackage{amsthm}
\usepackage{hyperref}
\usepackage[all]{xy}
\usepackage{todonotes}

\usepackage[english,catalan]{babel}

\theoremstyle{plain}
\newtheorem{teo}{Theorem}[section]
\newtheorem{lemma}[teo]{Lemma}
\newtheorem{lem}[teo]{Lemma}
\newtheorem{cor}[teo]{Corollary}
\newtheorem{prop}[teo]{Proposition}

\newtheorem{propes}[teo]{Properties}

\theoremstyle{definition}
\newtheorem{defi}[teo]{Definition}

\theoremstyle{remark}
\newtheorem{remark}[teo]{Remark}

\newenvironment{Proof}{\par\noindent{\bf Proof.}\enspace\ignorespaces}{\qed\par\par}

\newenvironment{pf1}{\par\noindent{\bf Proof of Theorem \ref{fundamental}}\enspace\ignorespaces}{\qed\par\par}
\newenvironment{pf2}{\par\noindent{\bf Proof of Corollary \ref{GraphtoGroup}}\enspace\ignorespaces}{\qed\par\par}

\let\emph\relax
\DeclareTextFontCommand{\emph}{\bfseries}

\newcommand{\bR}{{\mathbb{R}}}
\newcommand{\bP}{{\mathbb{P}}}
\newcommand{\bZ}{{\mathbb{Z}}}

\newcommand{\bC}{{\mathbb{C}}}

\newcommand{\bT}{{\mathbb{T}}}

\newcommand{\cO}{{\mathcal{O}}}
\newcommand{\cL}{{\mathcal{L}}}

\newcommand{\cT}{{\mathcal{T}}}

\newcommand{\mm}{\mathfrak{m}}
\newcommand{\BK}{\cT_K}
\newcommand{\AK}{\Aut(\bP^1_K)}
\newcommand{\AO}{\Aut(\bP^1_{\cO})}
\newcommand{\wK}{\bP^1(K)^{<3>}}
\newcommand{\Tog}{\Lambda}
\newcommand{\Bc}{B^{c}}
\newcommand{\ta}{\varpi}
\newcommand{\vp}{\gamma}
\newcommand{\mq}{\mu_q}

\newcommand{\vpB}[2]{\!\!\ ^{#1}{#2}}

\newcommand{\V}{V}
\newcommand{\E}{E}
\newcommand{\T}{T}

\DeclareMathOperator{\GL}{GL}\DeclareMathOperator{\Ends}{Ends}
\DeclareMathOperator{\PGL}{PGL} \DeclareMathOperator{\Aut}{Aut}
\DeclareMathOperator{\Star}{Star}\DeclareMathOperator{\Fix}{Fix}
\DeclareMathOperator{\id}{id} \DeclareMathOperator{\rank}{rank}
\DeclareMathOperator{\tr}{tr}

\begin{document}
\selectlanguage{english}

\author[Dani Samaniego]{Dani Samaniego}
\address{} \email{}

\author[Xavier Xarles]{Xavier Xarles}
\address{Departament de Matem\`atiques\\Universitat Aut\`onoma de
Barcelona\\08193 Bellaterra, Barcelona, Catalonia}
\email{xarles@mat.uab.cat}



\title{Schottky groups over complete valued fields and Bruhat-Tits $\Lambda$-trees.}

\begin{abstract}
Given a non-trivial complete valued field $K$ with value group
$\Tog$, we construct a $\Tog$-tree space associated to $K$ analog of
the Bruhat-Tits tree, and locally finite trees associated to compact
subset of $\bP^1(K)$. We propose a definition of hyperbolic matrix
and Schottky group over such field $K$. To any such Schottky group
$\Gamma$, we associate a compact set an action of $\Gamma$, such
that the quotient graph of the associated tree is a finite graph,
and $\Gamma$ is identified with its fundamental group. This results
extend the classical ones for discrete valuations of Mumford
\cite{mumford1972analytic} and non-archimedean rank 1 valuations of
Gerritzen and Van der Put \cite{GvdP}.
\end{abstract}

\maketitle \tableofcontents

In 1972, David Mumford in his celebrated paper
\cite{mumford1972analytic} constructed algebraic curves from certain
subgroups of $\PGL_2(K)$, for $K$ the field of fractions of a
complete noetherian local ring $\cO$. This construction imitates the
classical construction of Schottky over $\bC$, so the groups he
considered were called Schottky groups. Some years later, in 1980,
Lothar Gerritzen and Marius van der Put in his seminal book
\cite{GvdP} redid this construction but for fields complete with
respect to a non-archimedean absolute value. In this paper we
consider Schottky groups but now for a field $K$ complete with
respect to any valuation (which we may and do assume to be
non-trivial and written multiplicatively). The results we prove
generalize the ones obtained in the first chapter of \cite{GvdP}.
However we use a different construction, inspired by Vladimir
Berkovich analytic geometric \cite{berkovich2012spectral} introduced
in 1990, and concretely by the $\bR$-tree structure of the
analytification of $\bP^1_K$. See also \cite{GVX} for some results
on Mumford curves using these language. In our case we obtain a
$\Lambda$-tree, which is a generalization of a tree and of a
$\bR$-tree to a general totally ordered group $\Lambda$ (see
\cite{chiswell} for an introduction to $\Tog$-trees). We show that
this tree is the analogue of the classical Bruhat-Tits tree for
$\PGL_2(K)$ for $K$ discretely valued as introduced by Fran\c{c}ois
Bruhat and Jacques Tits \cite{bruhat1972groupes} and Jean-Pierre
Serre \cite{Serretrees}.

The structure of the paper is as follows: in the first section we
recall some well known results on valuations, and in the second we
introduce the $\Lambda$-tree of balls, which will substitute the
Berkovich tree in our general case. In the next section it is shown
that this tree is isomorphic to the (natural generalization of the)
Bruhat-Tits tree, obtained as a set as $\PGL_2(\cO)\backslash
\PGL_2(K)$, where $\cO$ denotes the associated valuation ring of $K$
(see \cite{MoSh}). We show after that how to associate a tree
$T(\cL)$ (and not just a sub-$\Lambda$-tree) to a compact subset
$\cL$ of $\bP^1(K)$, which we prove to be locally finite. The same
construction, but in the rank 1 case, was already done in
\cite{GvdP},(2.6.3). We also study how to recover the compact set
$\cL$ from the tree $T(\cL)$, in case $\cL$ is perfect, via the
classical theory of the ends of a tree. In the fourth section we
introduce and study what we call hyperbolic matrices in $\PGL_2(K)$.
The main difference with the classical case is that we insist some
defining element to be topologically nilpotent, which in the
classical rank 1 case is equivalent to have absolute value (in our
notation, valuation) less than 1, but in general is not. All these
is combined in section 5 to a definition of Schottky groups $\Gamma$
and its associated perfect and compact $\Gamma$-set $\cL_{\Gamma}$.
In section 6, which can be considered the core of the paper, we show
that the quotient of the tree $T(\cL_{\Gamma})$ with respect to the
natural action of $\Gamma$ is a finite graph, and that $\Gamma$ is
naturally identified with its fundamental group. Finally the last
section is devoted to adapt the construction given by Schottky to
our case, following the ideas in \cite{GvdP}, which also shows in
particular there are plenty of non-trivial such groups for any field
$K$ complete with respect to a non-trivial valuation. We end the
paper giving some applications of this construction. We tried to
write the paper as self contained as possible, thus reproving some
results which are probably well known, but for which we did not
found any clear reference. Some of the proofs are directly inspired
by the proofs of similar results in \cite{GvdP}; for others,
however, we tried to find more direct or more clear proofs, even for
the case considered there.

Comparing with the already mentioned first chapter of \cite{GvdP} we
do not consider any analogous concept of what they call
discontinuous groups, and hence we do not prove any result
concerning the existence of a normal Schottky subgroup with finite
index for a general finitely generated discontinuous group. This
type of results should be not difficult to obtain, and they can be
of independent interest in order to find Schottky groups in nature,
like for the case of the ones associated to Shimura or Drinfeld
modular curves.

\section*{Aknowledgements}
Part of this paper was the master thesis of the first author in
order to obtain the masters degree in Mathematics, under the
supervision of the second author and N\'{u}ria Vila. We thank
N\'{u}ria Vila for accepting to be the tutor of the first author,
and Angela Arenas for several comments and corrections.

The second author was partially supported by Generalitat de
Catalunya grant 2014-SGR-206 and Spanish MINECO grant
MTM2013-40680-P and MTM2016-75980-P.

\section*{Preliminaries and notations}\label{sect1}

In this section we collect some basic facts on valuations and
valuation fields. A general reference can be the book
\cite{EnglerPrestel}, but we will use here instead the
multiplicative notation.

Recall that a totally ordered (abelian) group with 0, $\Tog_0$, is
an abelian group $\Tog$ (which we will denote multiplicatively)
together with an absorbent element $0\notin \Tog$ verifying $0\cdot
\rho=\rho\cdot 0=0$ for all $\rho\in \Tog$, and with a total order
$\le$ such that
\begin{enumerate}
\item if $a\le b$ and $c\in \Tog$, then $a\cdot c\le b\cdot c$.
\item $0\le a$ for all $a\in \Tog$.
\end{enumerate}
We say that $\Tog_0$ is non-trivial if there exists $1\ne\rho\in
\Tog$.

A progression $\rho_n\in \Tog$ for $n\ge 1$ has limit $0$ if, for
every $\epsilon\in \Tog$, there exists $n_0\ge 1$ such that
$\rho_n\le \epsilon$ for all $n\ge n_0$. An element $\rho\in \Tog$
is topologically nilpotent (or a microbe) if the progression
$\rho^n$ has limit $0$ (the term microbe comes from the theory of
ordered fields, see e.g. \cite{brum}, page 89).

Given two elements $\rho_1\le \rho_2\in \Tog$, we denote the
intervals as usual
$$[\rho_1,\rho_2]:=\{\delta \in \Tog \ | \ \rho_1 \le \delta \le
\rho_2\}$$ and
$(\rho_1,\rho_2):=[\rho_1,\rho_2]\setminus\{\rho_1,\rho_2\}$. We
denote also $[\rho_1,\infty):=\{\delta \in \Tog \ | \ \rho_1 \le
\delta\}$.

\begin{defi}
A \textbf{surjective} map $|\ |:K\rightarrow\Tog\cup\{0\}$ from a
field $K$ to $\Tog\cup\{0\}$, where $\Tog$ is a non trivial ordered
group, is called a \emph{valuation} of $K$ if it satisfies
\begin{itemize}
\item $|xy|=|x|\cdot|y|$ $\forall x,y\in K$
\item $|x+y|\leq \max\{|x|,|y|\}$ $\forall x,y\in K$
\item $|x|=0\Leftrightarrow x=0$.
\end{itemize}
\end{defi}

Recall that if $|.|$ is a valuation of a field $K$ and if $a,b\in K$
are such that $|a|\neq |b|$ then $|a+b|=\max\{|a|,|b|\}$. Moreover
it is easy to see that $|1|=1$.

We say that a valuation $|.|$ has rank 1 if $\Tog_0 \hookrightarrow
\mathbb{R}_{\geq 0}$ as ordered groups. By composing with $-\log$ we
get the usual notion of (additive) real valuation.

Given a field $K$ with a valuation $|.|$, the ring of integers of
$K$ with respect to $|.|$ is $\cO_{|.|}=\{x\in K \ | \ |x|\leq 1\}$.
Note that $\cO_{|.|}$ is a local domain whose field of fractions is
$K$, with maximal ideal $\mm_{|.|}:=\{x\in K \ | \ |x|< 1\}$, and
residue field $k_{|.|}$. Moreover, the ring $\cO_{|.|}$ is a
valuation ring, and conversely any valuation ring $R$ is the ring of
integers of its field of fractions $K$ with valuation given by the
map projection map to $(K^*/R^*)\cup \{0\} $. For ease of
simplicity, we will denote it by $\cO$ if there is no risk of
confusion.

Given a sequence $\{a_n\}_n$ of elements in $K$, and $a\in K$, one
says that $\lim_{n\rightarrow\infty}a_n=a$ if and only if
$|a_n-a|\rightarrow 0$. It forms a Cauchy sequence if
$|a_{n+1}-a_n|\rightarrow 0$. Note that the notation can be
misleading as these notions depend on the given valuation.

The field $K$ is complete respect $|.|$ if every Cauchy sequence has
limit. Any field with a valuation can be subsumed into a minimal
field complete with respect to a valuation extending the given one,
called its completion (see for example \cite{EnglerPrestel}, Theorem
2.4.3.). Recall that any finite extension $L$ of a field $K$
complete with respect to a valuation has at least a valuation on it
extending the one of $K$, and moreover it is also complete.

We say that $q\in K^*$ is topologically nilpotent or a microbe (with
respect to $|.|$) if $|q|\in \Tog$ it is, i.e. if
$\lim_{n\rightarrow\infty} q^n=0$. If the valuation is of rank 1,
then $q$ is topologically nilpotent if and only if $|q|<1$; this is
never true if the rank is not 1. We say that the valuation is
microbial it there exists a microbe $q\in K^*$ (see \cite{HuEt},
definition 1.1.4.).


A valuation ring $\cO$ and his field of fractions $K$ inherits a
natural topology which make them topological rings. In case $K$ has
a microbe $q$, it can be described as the $I$-adic topology for
$I=q\cO$, which is independent of $q$. A basis of open sets is
formed by the (closed) balls with radius $\rho\in\Tog$.

Given $p\in K$ and $\rho\in\Tog_0$, the (closed) ball with center
$p$ and radius $\rho$ is $B(p,\rho)=\{y\in K \ | \ |y-p|\leq
\rho\}$.

When considering the projective line $\bP^1(K)=K\cup\{\infty\}$ with
its inherited (analytic) topology, the set of closed balls is not a
basis for the topology; one needs to include also the complements of
the open balls
$$\Bc(p,\rho):=\{z\in K \ | \ |z-p|\geq\rho\}\cup\{\infty\}$$
for $p\in K$ and $\rho\in\Tog$ to get a subbasis. Given $p\in K$,
$\rho_1\in \Tog_0$ and $\rho_2\in \Tog\cup\{\infty\}$, the
(generalized) annulus with center $p$ and radii $\rho_1$ and
$\rho_2$ is  $$C(p,\rho_1,\rho_2):=\{z\in \bP^1(K) \ | \
\rho_1\leq|z-p|\leq\rho_2\}.$$ Note that the case $\rho_1=0$ are
closed balls, and $\rho_2=\infty$ are complements of open balls. The
set of all generalized annulus form a basis for the topology of
$\bP^1(K)$. This is proven by observing that the intersection of two
generalized annulus is either empty or a generalized annulus.

We will denote also by $C(p,\rho):=C(p,\rho,\rho)=B(p,\rho)\cap
\Bc(p,\rho)$ the circle with center $p$ and radius $\rho$.

\section{The $\Lambda$-tree of Balls}\label{sect2}

\begin{defi}
We define the \emph{space of balls} $\BK=\{B(p,\rho) \ | \ p\in K,
\rho\in \Tog \}$. We also define $\overline{\BK}=\BK\cup
K\cup\infty$, which can be seen as the set of balls with radius
$\rho\in \Tog_0\cup\{\infty\}$, being $B(p,\infty):=K$ for all $p\in
K$.
\end{defi}

We will see in this section that $\BK$ has a natural structure of
(oriented) $\Lambda$-tree, a generalization of (simplicial) trees
and $\bR$-trees. The order will be the inclusion relation (as a
subsets of $K$). The main property is given by the following
elementary result.

\begin{lemma} Define $\varrho:\BK\to \Tog$ by
$\varrho(B(p,\rho)):=\rho$. Then $\varrho$ is well-defined and for
any $B(p,\rho)\in \BK$ induces a bijection
$$\{B\in \BK\ | \ B\subset B(p,\rho)\}\cong [\rho,\infty).$$
\end{lemma}

\begin{Proof} We know that $B(x_1,\rho_1)=
B(x_2,\rho_2)$ if and only if $\rho_1=\rho_2$ and $x_2\in
B(x_1,\rho_1)$, which shows that $\varrho$ is well defined, and that
for any $\delta\ge \varrho$, $$\{B\in \BK\ | \ B\subset
B(p,\rho)\}\cap \varrho^{-1}(\delta)=\{B(p,\delta)\}.$$
\end{Proof}

\begin{lemma} Given two balls $B_1$ and $B_2\in \BK$, there exists a
minimal ball $B_1\vee B_2$ that contains both. Even more, the set
$$\{ B\in \BK \ | \ B\supset B_1 \mbox{ and } B\supset B_2\}$$
is totally ordered and with a minimal element with respect to the
inclusion.
\end{lemma}

\begin{Proof} We only need to observe that
$$B(x_1,\rho_1)\vee
B(x_2,\rho_2)=B(x_1,\max\{|x_1-x_2|,\rho_1,\rho_2\})$$ verifies the
properties. The second assertion is due to the well-known property
asserting that $B(x_1,\rho_1)\cap B(x_2,\rho_2)$ is either empty or
the smallest of both.
\end{Proof}

We can give now a structure of $\Tog$-metric space for $\BK$.
\begin{defi}\label{distance} Given two balls $B_1\subset B_2$ of $\BK$, we define
the $\Tog$-distance between them as
$d(B_1,B_2):=\varrho(B_2)\varrho(B_1)^{-1}$. Given any two balls
$B_1$, and $B_2$ in $\BK$, we define
$$d(B_1,B_2):=d(B_1,B_1\vee B_2)d(B_2,B_1\vee B_2).$$
\end{defi}

The following properties are elementary, and they define and show
that $\BK$ is a $\Tog$-metric space.

\begin{propes} The map $d\colon \BK\times \BK \to \Tog$ verifies
\begin{enumerate}
\item $d(B_1,B_2)\ge 1$ for any $B_1$ and $B_2\in \BK$.
\item $d(B_1,B_2)=1$ if and only if $B_1=B_2$.
\item $d(B_1,B_2)=d(B_2,B_1)$  for any $B_1$ and $B_2\in \BK$.
\item $d(B_1,B_2)\le d(B_1,B_3)d(B_3,B_2)$ for any $B_1$, $B_2$ and $B_3 \in \BK$.
\end{enumerate}
\end{propes}

Note that we are considering the not so usual (in this context)
multiplicative notation. We can then define segments in $\BK$ and
show it is geodesically linear: given any two balls there is a
unique segment going from one to the other. This shows that $\BK$ is
a $\Tog$-tree as defined in \cite{chiswell}.

\begin{remark} In the case $|.|$ is a non-archimedean valuation, the
space $\BK$ is form by the so called type II points inside
$\bP^{1,an}_K$, with its natural metric (see
\cite{berkovich2012spectral}). This form all the points of
$\bP^{1,an}_K\setminus \bP^1(K)$ if and only if $|.|$ has image all
$\bR_{\ge 0}$ and $K$ is spherically complete.
\end{remark}

Recall that a segment $[B_1,B_2]_{\Tog}$ from $B_1$ to $B_2$ is an
isometry $\alpha\colon [\rho_1,\rho_2]\to \BK$ where $\rho_1$ and
$\rho_2\in \Tog$ such that $\alpha(\rho_1)=B_1$ and
$\alpha(\rho_2)=B_2$.

We define the path from $p\in K$ to $\infty$ inside $\BK$ as
$\pi(p,\infty)=\{B(p,\rho) \ | \ \rho\in \Tog \}.$ In general, given
any ball $B(p,\rho)$ for $\rho\in \Tog\cup\{\infty\}$, we define the
path
$$\pi(p,B(p,\rho))=\{B(p,\delta) \ | \ \delta\le \rho\}.$$

Note that $\varrho$ induces an isometry $\varrho\colon
\pi(p,\infty)\cong \Tog$. The intersection of two such paths is
clearly $$
    \pi(p,\infty)\cap\pi(q,\infty)=\{B(p,r) \ | \ r\geq |p-q|\}=:\pi(B(p,|p-q|),\infty))\cong[|p-q|,\infty)
$$ and we define
$$
    \pi(p,q):=\pi(p,B(p,|p-q|))\cup\pi(q,B(p,|p-q|)).
$$

\begin{lemma}
Let $p_1$, $p_2$ and $p_3$ three distinct points in $\bP^1(K)$. Then
$$\pi(p_1,p_2)\cap\pi(p_2,p_3)\cap\pi(p_1,p_3)=\{t(p_1,p_2,p_3)\}$$
a unique point $t(p_1,p_2,p_3)$ of $\BK$. Note that
$t(p_1,p_2,\infty)=B(p_1,|p_1-p_2|)$. In general, if
$|p_1-p_3|=|p_2-p_3|\geq |p_1-p_2|$ then
$t(p_1,p_2,p_3)=B(p_1,|p_1-p_2|)$.
\end{lemma}

\begin{Proof}
First observe that the case one of the points is $\infty$, say
$p_3=\infty$, is clear from the definition. Second, if all points
are in $K$, we can suppose that $|p_1-p_3|=|p_2-p_3|\geq |p_1-p_2|$
after changing the order if necessarily. Then $B(p_1,|p_1-p_2|)\in
\pi(p_1,p_2)\cap\pi(p_2,p_3)\cap\pi(p_1,p_3)$ clearly. On the other
hand, a ball $B(p_1,\delta)\in \pi(p_1,p_3)$ if $\delta\leq
|p_1-p_3|$, and then $B(p_2,\delta)\in \pi(p_2,p_3)$ since
$\delta\leq |p_2-p_3|=|p_1-p_3|$. But then
$B(p_2,\delta)=B(p_1,\delta)$ if $\delta\geq |p_1-p_2|$, and
$B(p_1,\delta)\in \pi(p_1,p_2)$ if $\delta\leq |p_1-p_2|$; so
$\delta=|p_1-p_2|$ is the unique solution.
\end{Proof}

From the proof of the lemma we can see that for any three distinct
points $p_1$, $p_2$ and $p_3$ in $K$, if we order them such that
$|p_1-p_3|=|p_2-p_3|\geq |p_1-p_2|$, then
$t(p_1,p_2,p_3)=t(p_1,p_2,\infty)$. We call such ordering a ball
ordering.

Note that we have a natural bijection $i_{p_1,p_2}$ from
$\pi(p_1,p_2)$ to $\Tog$, once fix an ordering of $p_1$ and $p_2$,
and which sends the ball $B(p,|p-q|)$ to 1, given by
$i_{p_1,p_2}(t(p_1,p_2,q))=\frac{|q-p_1|}{|q-p_2|}$.

We denote $$\wK:=\{(p_1,p_2,p_3)\in \bP^1(K)^3 / p_1\neq p_2\neq
p_3\neq p_1\}=(\bP^1(K)^3\backslash\Delta),$$ where the set $\Delta$
is formed by the points $(p_1,p_2,p_3)$ such that  $p_1\neq p_2\neq
p_3\neq p_1$.

So we have $t:\wK \to \BK$ as $t(p_1,p_2,p_3):=B(p_i,\rho)$ where
$\rho$ is the smallest distance between the three points and $p_i$
is one of the two points that gives this smallest distance.

\section{The tree of balls and the Bruhat-Tits tree}

In this section we show that the $\Tog$-tree constructed in the
previous section coincides with the Bruhat-Tits tree, constructed in
the language of lattices in \cite{MoSh}, II.3..

Consider the group of automorphims $\AK\cong \PGL_2(K)$ of the
projective line over $K$, which we will identify with the projective
linear group of matrices through the usual isomorphism. Recall that
there is a natural bijection between $\wK:=\{(p_1,p_2,p_3)\in
\bP^1(K)^3 / p_1\neq p_2\neq p_3\neq p_1\}$ and $\AK$, given by
sending $\varphi\in\AK$ to the triple
$(\varphi(0),\varphi(1),\varphi(\infty))$. The group $\AK$ acts on
$\wK$ by the usual left action $\phi(\varphi)=\phi\circ \varphi$,
which is clearly transitive. Explicitly it is given by, for any
$\tau\in \AK$,
$\tau(t(p_1,p_2,p_3)):=t(\tau(p_1),\tau(p_2),\tau(p_3))$. In this
section we will show that this action descents to an action of $\AK$
on $\BK$ via the $t$-map, and that this action gives an
identification of $\BK$ as an analogous of the Bruhat-Tits
($\Lambda$-)tree of $K$ (with respect to the valuation).

To understand how the action descends it is natural to give an
alternative description using balls. But it is clearly not true that
$\tau(B)$ is a ball for any ball $B$ and for any $\tau\in \AK$, as
it shows the example of $\tau(t)=1/t$ and $B=B(0,1)$, since
$\infty\in \tau(B)$. The following shows that this is in fact the
only obstruction.

If $\gamma\in \AK$, we denote by $\gamma'(p):=(bc-ad)(cp+d)^{-2}$
the derivative of $\gamma(t):=(at+b)/(ct+d)$ with respect to $t$
applied to $p$. Denote also $$\gamma'(\infty):=(bc-ad)c^{-2}.$$

\begin{lem}\label{taylor} Consider $\gamma \in \AK$, $\gamma\ne \id$, $p\in K$ and
$\delta\in \Tog$. Suppose that $\infty\notin \gamma(B(p,\delta))$.
Then
$$|\gamma(p)-\gamma(q)|=|\gamma'(p)||p-q|$$ for all $q\in
B(p,\delta)$.
\end{lem}

\begin{Proof} First, note that, if
$\gamma(t)=(at+b)/(ct+d)$, then
$$\gamma(p)-\gamma(q)=(p-q)\gamma'(p)\frac{cp+d}{cq+d}.$$
Now, the condition  $\infty\notin \gamma(B(p,\delta))$ is equivalent
to $|cp+d|>\delta |c|$. But then $|cp+d|=|cq+d|$ since
$|(cp+d)-(cq+d)|=|c||p-q|\le |c|\delta <|cp+d|$.
\end{Proof}

\begin{cor}\label{actiononballs} Consider $\gamma \in \AK$, $\gamma\ne \id$, $p\in K$ and
$\delta\in \Tog$.
\begin{enumerate}
\item Suppose that $\infty\notin
\gamma(B(p,\delta))$. Then $\gamma(B(p,\delta))=B(\gamma(p),
|\gamma'(p)|\delta)$.
\item Suppose that $\infty\in
\gamma(B(p,\delta))$. Then
$\gamma(B(p,\delta))=\Bc(\gamma(\infty),|\gamma'(\infty)|\delta^{-1})$.
\end{enumerate}
\end{cor}

\begin{Proof} First of all, note that, if $\gamma(t)=(at+b)/(ct+d)$,
then $\infty\notin \gamma(B(p,\delta))$ is equivalent to
$|cp+d|>\delta |c|$.

Suppose first that  $\infty\notin \gamma(B(p,\delta))$. If $q\in
B(p,\delta)$, then, by lemma \ref{taylor},
$$|\gamma(p)-\gamma(q)|=|\gamma'(p)||p-q|\le|\gamma'(p)|\delta$$
which shows that $\gamma(q)\in B(\gamma(p), |\gamma'(p)|\delta)$. On
the other side, we first show that $\infty\notin
\gamma^{-1}(B(\gamma(p), |\gamma'(p)|\delta))$, or, equivalently,
that $a/d=\gamma(\infty)\notin B(\gamma(p), |\gamma'(p)|\delta)$.
But
$$|\gamma(p)-\gamma(\infty)|=|\gamma'(p)|\frac{|cp+d|}{|c|}>\delta.$$
Hence we can apply again lemma \ref{taylor} and we get that for any
$q'\in B(\gamma(p), |\gamma'(p)|\delta)$,
$$|\gamma^{-1}(q')-p|=|(\gamma^{-1})'(\gamma(p))||q'-\gamma(p)|\le|(\gamma^{-1})'(\gamma(p))|
|\gamma'(p)|\delta=\delta$$ by chain's rule, which shows 1.

The second assertion is shown with analogous arguments. If
$\infty\in \gamma(B(p,\delta))$, then
$$|\gamma(p)-\gamma(\infty)|=|\gamma'(\infty)|\frac{|c|}{|cp+d|}\le |\gamma'(\infty)|\delta^{-1}$$
and the same works for any $q\in B(p,\delta)$, so
$\gamma(B(p,\delta))\subset
\Bc(\gamma(\infty),|\gamma'(\infty)|\delta^{-1})$. Now
$\gamma^{-1}(\infty)\in B(p,\delta)$ since
$$|\gamma^{-1}(\infty)-p|=\frac{|cp+d|}{|c|}\le \delta$$
by hypothesis. Hence $B(p,\delta)=B(\gamma^{-1}(\infty),p)$. But for
any $q\in\Bc(\gamma(\infty),|\gamma'(\infty)|\delta^{-1})$, so
$|cq-a|\ge |c||\gamma'(\infty)|\delta^{-1}$, we get that
$$|\gamma^{-1}(q)-\gamma^{-1}(\infty)|=|(\gamma^{-1})'(\infty)|\frac{|c|}{|cq-a|}\le
\frac{|(\gamma^{-1})'(\infty)|}{|\gamma'(\infty)|\delta^{-1}}=\delta$$
which shows the reverse inclusion.
 \end{Proof}

\begin{prop}\label{Actionont} For any $(p_1,p_2,p_3)\in\wK$
and for any $\gamma\in\AK$, let $B=B(p,\delta)$ be the ball
representing $t(p_1,p_2,p_3)$.
\begin{enumerate}
\item If $\infty\notin \gamma(B)$, then $\gamma(B)$ represents
$t(\gamma(p_1),\gamma(p_2),\gamma(p_3))$.
\item If $\infty\in \gamma(B)$, then
$B(\gamma(\infty),|\gamma'(\infty)|\delta^{-1})$ represents
$t(\gamma(p_1),\gamma(p_2),\gamma(p_3))$.
\end{enumerate}
As a consequence the action of $\AK$ on $\wK$ descents to an action
on $\BK$.
\end{prop}
\begin{Proof} We can and will suppose $(p_1,p_2,p_3)$ are in ball position, so
$|p_1-p_2|\le |p_1-p_3|=|p_2-p_3|$, and $B=B(p_1,|p_1-p_2|)$. We
denote as above $\gamma(t)=(at+b)/(ct+d)$.

Suppose first $\infty\notin \gamma(B)$, so
$\gamma(B)=B(\gamma(p_1),|\gamma'(p_1)||p_1-p_2|)$. Also we have by
lemma \ref{taylor}
$|\gamma(p_1)-\gamma(p_2)|=|\gamma'(p_1)||p_1-p_2|$. But
$$|\gamma(p_1)-\gamma(p_3)|= |p_1-p_3||\gamma'(p_1)|\frac{|cp_1+d|}{|cp_3+d|}$$
is equal to  $|\gamma(p_2)-\gamma(p_3)|$, since
$|\gamma'(p_1)|=|\gamma'(p_2)|$ and $|cp_1+d|=|cp_2+d|$. Hence
$\gamma(p_1),\gamma(p_2),\gamma(p_3)$ are in ball position and
$B(\gamma(p_1),|\gamma(p_1)-\gamma(p_2)|)=\gamma(B)$.

Now, if $\infty\in \gamma(B)$, so
$\gamma(B)=\Bc(\gamma(\infty),|\gamma'(\infty)|\delta^{-1})$. Denote
$p':=\gamma(\infty)$ and $\delta':= |\gamma'(\infty)|\delta^{-1}$,
and we want to show that $B(p',\delta')$ represents
$t(\gamma(p_1),\gamma(p_2),\gamma(p_3))$. Note that we have
$$|\gamma(p_i)-p'|\ge \delta' $$
for $i=1,2$. We divide the proof in two cases.

Suppose first that $p_3\notin B(p_1,\delta)$, so
$|p_3-p_2|>|p_2-p_1|=\delta$. Then
$$|cp_3+d|=|c||p_3-p_1|>|c|\delta\ge |cp_1+d|$$
hence $|\gamma(p_3)-\gamma(\infty)|<\delta'$, thus $\gamma(p_3)\in
B(p',\delta')$. Now, observe that either $|cp_1+d|=|c|\delta$ or
$|cp_2+d|=|c|\delta$, since $|(cp_1+d)-(cp_2+d)|=|c|\delta$ and both
are $\le |c|\delta$. We can and will suppose it is $p_1$. But then
$$|\gamma(p_3)-\gamma(p_1)|=|\gamma'(\infty)|\frac{|c|}{|cp_1+d|}=\delta'$$
and $|\gamma(p_3)-\gamma(p_2)|\ge \delta'$, which shows that
$\gamma(p_3)$, $\gamma(p_1)$ and $\gamma(p_2)$ are in ball position
and $B(\gamma(p_3),|\gamma(p_3)-\gamma(p_1)|)=B(p',\delta')$.

Now, suppose that  $p_3\in B(p_1,\delta)$, hence
$|p_3-p_2|=|p_1-p_2|=|p_3-p_1|=\delta$. Same arguments as in the
previous case show that there exists $i$ such that $|cp_i+d|\le
\delta|c|$ and $|cp_j+d|\le \delta|c|$ for $j\ne i$. We can and will
suppose that $i=1$. Then one shows that
$|\gamma(p_3)-\gamma(\infty)|=\delta'$, and the same for $p_2$, that
$|\gamma(p_3)-\gamma(p_2)|=\delta'$, and finally that
$|\gamma(p_3)-\gamma(p_1)|=|\gamma(p_2)-\gamma(p_1)|\ge \delta'$,
which implies the result.
\end{Proof}

Given a ball $B$, with corresponding $t\in \BK$, and an automorphism
$\vp$, we will denote $\vpB{\vp}{B}$ the ball corresponding to
$\vp(t)$. So $\vpB{\vp}{B}=\vp(B)$ if this last set is a ball, or
$\vpB{\vp}{B}=B(\vp(\infty),|\vp'(\infty)|\varrho(B)^{-1})$ if it is
not. The automorphisms of $\AK$ preserve the distance between balls,
using definition \ref{distance}, hence they are isometries.

\begin{lem}\label{pairofballs} For any pair of balls $B$ and $B'$, and
an automorphism $\vp\in \AK$, we have
$d(\vpB{\vp}{B},\vpB{\vp}{B'})=d(B,B')$.
\end{lem}

\begin{Proof} Suppose first that $B'=B(p,\delta')\subset B=B(p,\delta)$.
If $\vp(B)$ is a ball, then we have $\vp(B')\subset \vp(B)$.
Corollary \ref{actiononballs} asserts that
$\varrho(\vp(B))=|\vp'(p)|\delta$ and
$\varrho(\vp(B'))=|\vp'(p)|\delta'$. Then
$$d(\vp(B'),\vp(B))=\varrho(\vp(B))\varrho(\vp(B'))^{-1}=\frac{|\vp'(p)|\delta}{|\vp'(p)|\delta'}=\frac{\delta}{\delta'}=d(B,B').$$
If $\vp(B')$ is not a ball, applying corollary \ref{actiononballs}
we get that
$\vpB{\vp}{B'}=B(\vp(\infty),|\vp'(\infty)|{\delta'}^{-1})\supset
\vpB{\vp}{B}=B(\vp(\infty),|\vp'(\infty)|{\delta}^{-1})$, hence
$d(\vpB{\vp}{B},\vpB{\vp}{B'})=\delta/\delta'=d(B,B')$ as above.

If $\vp(B)$ is not a ball, but $\vp(B')$ it is, then
$\vpB{\vp}{B}=B(\vp(\infty),\delta')$, where
$\delta'=|(\vp)'(\infty)|\delta^{-1}$ by corollary
\ref{actiononballs}. Now, $\vp(B')\vee
\vpB{\vp}{B}=B(\vp(p),|\vp(p)-\vp(\infty)|)$, and
$$d(\vp(B'),\vpB{\vp}{B})=d(B(\vp(p),|\vp(p)-\vp(\infty)|),\vp(B'))d(B(\vp(p),|\vp(p)-\vp(\infty)|),\vpB{\vp}{B})=$$
$$=\frac{|\vp(p)-\vp(\infty)|^2}{|\vp'(p)||(\vp)'(\infty)|}\frac{\delta}{\delta'}.$$
But one easily shows that
$|\vp(p)-\vp(\infty)|^2=|\vp'(p)||(\vp)'(\infty)|$.
\end{Proof}

We denote by
$\ta(\varphi)=t(\varphi(0),\varphi(1),\varphi(\infty))$. The
equivalence relation determined by $\varphi\sim\varphi'$ when
$\ta(\varphi)=\ta(\varphi')$ is thus determined by the stabilizer of
an element, say $t_0:=t(0,1,\infty)=\cO$.

\begin{teo} An automorphism $\varphi\in \AK$ stabilizes $t_0$ if and
only if $\varphi\in \AO$.
\end{teo}

\begin{Proof} Denote by $\Gamma_0\subset\Aut(\bP^1_{K})$ the stabilizer
of $t_0$. Observe that an automorphism $\varphi\in \Aut(\bP^1_{K})$
is in fact $\varphi\in \AO$ if and only if it can be written as
$\varphi(t)=(at+b)/(ct+d)$ for $a,b,c,d\in \cO$ with $|ad-bc|=1$.

First of all, observe that the automorphisms $\tau$ that fix the set
$\{0,1,\infty\}$ are both in $\Gamma_0$ and also in $\AO$. If we
compose an automorphism $\psi$ with one such $\tau$ in order to
obtain an automorphism $\vp=\tau\circ\psi$ that $\vp(0)$, $\vp(1)$
and $\vp(\infty)$ are ball ordered, so
$t(\vp(0),\vp(1),\vp(\infty))=t(\vp(0),\vp(1),\infty)=B(\vp(0),|\vp(0)-\vp(1)|)$,
then $\psi\in\Gamma_0$ (respectively $\psi\in \AO$) if and only if
$\vp\in \Gamma_0$ (respectively $\vp\in \AO$). So we are reduced to
consider only the case that $\vp$ verifies that $\vp(0)$, $\vp(1)$
and $\vp(\infty)$ are ball ordered, which we will say that $\vp$ is
ball suited.

We will show first that $\Gamma_0\subset \AO$. We decompose a ball
suited $\vp$ as composition of two automorphisms: the automorphism
$\vp_0$ which sends $\vp_0(0)=\vp(0)$, $\vp_0(1)=\vp(1)$ and
$\vp_0(\infty)=\infty$, and the diagonalizable automorphism $\vp_1$,
sending $\vp_1(\vp(0))=\vp(0)$, $\vp_1(\vp(1))=\vp(1)$ and
$\vp_1(\infty)=\vp(\infty)$. Hence $\vp_1$ has two fixed points,
$\vp(0)$ and $\vp(1)$. We will see that if $\vp_i\in \Gamma_0$ then
$\vp_i\in \AO$ for $i=1,2$, and that $\vp\in \Gamma_0$ if and only
if $\vp_i\in \Gamma_0$ for $i=0$ and $1$. Note that $\Gamma_0$ and
$\AO$ are subgroups, hence, if the $\vp_i$ are in one of them for
both $i=0$ and $1$, so it is their composition.

Now, if $\vp\in\Gamma_0$, then
$\vp_0(t_0)=t(\vp(0),\vp(1),\infty)=t(\vp(0),\vp(1),\vp(\infty))=t_0$,
being $\vp$ ball suited. Hence $\vp_0\in \Gamma_0$, and thus
$\vp_1\in \Gamma_0$.

Therefore we are reduced to show the result for the automorphisms of
the type $\vp_0$ and $\vp_1$. The first case is easy; one has
$\vp_0(t)=at+b$, with $a=\vp_0(1)-\vp_0(0)$ and $b=\vp(0)$. Since
$\vp_0$ is ball suited, $t(\vp(0),\vp(1),\infty)=t_0$ if and only if
$|a|=1$ and $|b|\le 1$, which happens exactly when $\vp_0\in\AO$.

We consider now the second case of diagonalizable automorphisms of
the type $\vp_1$, with fixed points $p_0$ and $p_1$, and
 $\vp_1\in \Gamma_0$, hence $t(p_0,p_1,\infty)=t_0$. Therefore $|p_0|\le
1$, $|p_1|=1$ and $|p_0-p_1|=1$. Take $\tau\in \AK$ such that
$\tau(0)=p_0$, $\tau(\infty)=p_1$ and $\tau(1)=\infty$. Explicitly
$\tau(t)=(p_0t-p_1)/(t-1)$, which is clearly in $\AO$ since
$|p_1-p_0|=1$. Moreover $\tau\in\Gamma_0$ since
$\tau(t_0)=t(p_0,\infty,p_1)=t_0$. Then $\psi:=\tau^{-1}\vp_1\tau\in
\Gamma_0$ verifies $\psi(0)=0$ and $\psi(\infty)=\infty$, so
$\psi(t)=qt$ for some $q\in K^*$. But then $\psi\in\Gamma_0$ if and
only if $1=|\psi(1)|=|q|$. Hence $\psi\in\AO$, so
$\vp_1=\tau\psi\tau^{-1}$ also.

Finally, suppose $\vp(t)=(at+b)/(ct+d)\in \AO$ with $a,b,c$ and
$d\in \cO$ and $|ad-bc|=1$ and we want to show $\vp(t_0)=t_0$. If
$\infty\notin\vp(B(0,1))$, equivalently $|d|>|c|$, so $|d|=|a|=1$,
then $|\vp(0)|=|c/d|<1$ and $|\vp'(0)|=|d|^{-2}=1$. Therefore
$\vp(B(0,1))=B(\vp(0),1)=B(0,1)$ by corollary \ref{actiononballs}.
If, on the contrary, $\infty\in\vp(B(0,1))$, i.e. $|d|\le |c|$, so
$|c|=1$, then $\vp(\infty)=a/c\in\cO$ and
$|\vp'(\infty)|=|c|^{-2}=1$. Therefore, using again the same
corollary, $\vp(B(0,1))=\Bc(0,1)$. Proposition \ref{Actionont}
implies the result.
\end{Proof}

The following corollary is a well-known consequence of the
transitivity of the action of $\AK$ on $\BK$ and the previous
theorem.

\begin{cor} The map $\ta:\AK\to \BK$ determines a
bijection $$\AO\backslash \AK \cong \BK.$$
\end{cor}

\section{The tree associated to a compact set}



\begin{defi}
Let $\cL\subset\bP^1(K)$ be with at least three elements. We define
the $\Tog$-subtree associated to $\cL$ as
\begin{equation*}
    \bT(\cL)=\bigcup_{\substack{p_1,p_2\in\cL \\ p_1\neq p_2}}\pi(p_1,p_2)
\end{equation*}
\end{defi}

Given the $\Tog$-tree $\bT(\cL)$ we will construct a (simplicial)
graph as follows: the set of vertices is
\begin{equation*}
    \V(\cL)=\{t(p_1,p_2,p_3) \ | (p_1,p_2,p_3)\in\cL^3\cap\wK \}\subset \bT(\cL).
\end{equation*}
Given  $v_1,v_2\in \V(\cL)$, we say that they determine an edge
$[v_1,v_2]$  with ends $v_1$ and $v_2$ if
\begin{equation*}
    [v_1,v_2]_{\Tog}\cap \V(\cL)=\{v_1,v_2\}.
\end{equation*}
We denote the set of edges as $\E(\cL)$. The (simplicial) graph they
determine will be denote by $\T(\cL)$. It is easy to show it is a
forest (a union of trees).

Note that, if $\cL$ is finite, then the graph $\T(\cL)$ is a tree.
This is because, given any two vertices $v_1$ and $v_2$ of the
graph, the segment $[v_1,v_2]_{\Tog}$ can be subdivided in a finite
number of edges.

Note that, if $\cL$ contains exactly three points, then $\T(\cL)$
has only one vertex and no edge, whereas if it has four points, then
it has either one vertex or two vertices and one edge. In general
the number of vertices is bounded by $\# \cL-2$ and the number of
edges by $\# \cL-3$. The following lemma shows this result by
induction and it will be useful later.

\begin{lemma}\label{adding} Given $\cL\subset \bP^1(K)$ any subset, the
inclusion $\cL\subset \overline{\cL}$ into the closure gives natural
bijections $\bT(\cL)=\bT(\overline{\cL})$ and
$\T(\cL)=\T(\overline{\cL})$.

If now $\cL$ is closed, $\# \cL\ge 2$ and $p\in \bP^1(K)\setminus
\cL$, then $\V(\cL)\cup \{ v_p\}=\V(\cL\cup \{p\})$ for a vertex
$v_p$ (which may or may not be in $\V(\cL)$).
\end{lemma}

\begin{Proof} Given $p_1$ and $p_2\in \overline{\cL}$, which we will suppose
are not equal to $\infty$, and two balls $B_1$ and $B_2 \in
\pi(p_1,p_2)$, first we will show that $[B_1,B_2]_{\Tog}\subset
\bT(\cL)$. We can and will suppose that $B_1=B(p_1,\delta_1)$ and
$B_2=B(p_2,\delta_2)$, with $\delta_1$ and $\delta_2\le |p_1-p_2|$.
Then there exists $p_1'\in \cL\cap B(p_1,\delta_1)$ and $p_2'\in
\cL\cap B(p_2,\delta_2)$, and then $[B_1,B_2]_{\Tog}\subset
\pi(p_1',p_2')\subset \bT(\cL)$.

Now, suppose moreover than $B_1\in \V(\overline{\cL})$, so
$B_1=t(p_1,p_2,p_3)$ for some $p_1, p_2$ and $p_3\in
\overline{\cL}$. Taking as before three points $p_i'\in \cL$
sufficiently close to $p_i$ for all $i=1,2,3$, we have
$t(p_1,p_2,p_3)=t(p_1',p_2',p_3')\in \V(\cL)$.

If one of the points is equal to $\infty$, one adapts the argument
by using complements of open balls.

Finally, to show the last assertion, suppose $p\ne \infty$. Then
$v_p=B(p,\delta)$, where $$\delta:=\sup\{\epsilon\in \Tog\ |
B(p,\epsilon)\cap \cL =\emptyset\}.$$
\end{Proof}

We will show that $\T(\cL)$ is a tree for any compact subset.

\begin{teo}\label{teofinite} Let $\cL\in\bP^1(K)$ be a compact subset.
Then $[v_1,v_2]_{\Tog}\cap \V(\cL)$ is finite for any $v_1$ and
$v_2\in \V(\T(\cL))$ .
\end{teo}

\begin{Proof} We will suppose $\infty\in \cL$, since if it is not then
$\cL\cup\{\infty\}$ would be also compact, and $\V(\cL)\subset
\V(\cL\cup\{\infty \})$ (and it fact it contains at most one more
vertex). Given two vertices $v_1$ and $v_2$, we denote by $v_1\vee
v_2$ the element in $\BK$ corresponding to the minimal ball
containing both. Since $\infty\in \cL$, $v_1\vee v_2\in V(\cL)$.
Then
$$[v_1,v_2]_{\Tog}=[v_1,v_1\vee v_2]_{\Tog}\cup[v_1\vee v_2,v_2]_{\Tog}$$
hence we are reduced to show only the case that $v_1\le v_2$ with
respect to the partial order of $\BK$.

Then  $v_1=B(p,\rho_1)\subsetneq B(p,\rho_2)=v_2$ and
\begin{equation*}
    [v_1,v_2]_{\Tog}\cap \V(\T(\cL))=\{B(p,\rho) \ | \ \rho_1\leq\rho\leq\rho_2\}
\end{equation*}
where $B(p,\rho)=t(p,q,\infty)$ for some $q\in\cL$. We want to see
that there are a finite number of such $q$. Each $q$ is in the
circle $C(p,\rho)$ for $\rho_1\leq\rho\leq\rho_2$. We have
\begin{equation*}
 \cL\subset B(p,\rho_1)\cup \Bc(p,\rho_2) \cup \bigcup_{\rho_1<\rho<\rho_2}C(p,\rho)
\end{equation*}
with  $C(p,\rho)\cap\cL\neq\emptyset$. The sets are disjoint and
open, so there is a finite number of them.
\end{Proof}

\begin{cor}
The graph $\T(\cL)$ is a tree.
\end{cor}

\begin{Proof}
We need to show it is connected. But given two vertices $v$ and
$v'$, we have $$[v,v']_{\Tog}\cap
\V(\cL)=\{v=v_0,v_1,\cdots,v_n,v_{n+1}=v'\}=[v,v_1]\cup\cdots\cup[v_n.v']$$
for some $n\ge 0$, and any of these $[v_i,v_{i+1}]$ are edges.
Clearly this is the unique path from $v$ to $v'$, hence it is a
tree.
\end{Proof}

Recall that the star of a vertex $v\in \V(\cL)$ is
$\Star_{\T(\cL)}(v)=\{[v,w]\in \E(\cL)\}$.

A graph is called \emph{locally finite} if $\Star_{\T(\cL)}(v)$ is
finite for all $v\in \V(T)$.

\begin{teo} The tree $\T(\cL)$ is locally finite.
\end{teo}

\begin{Proof}
As in the proof of theorem \ref{teofinite}, we will construct a
covering of our compact set $\cL$ by non-empty and disjoint open
sets indexed by $\Star_{\T(\cL)}(v)=\{[v,w]\in \E(\cL)\}$, at least
when $\cL$ has no isolated points. Let $\cL^i$ be the set of
isolated points of $\cL$; since $\cL$ is compact, $\cL^i$ is finite.
Consider $\cL':=\cL\setminus \cL^i$, which is also compact. Then
Lemma \ref{adding} allows us to show that $\T(\cL)$ is locally
finite if and only if $\T(\cL')$ is locally finite. So we can and
will suppose that $\cL$ has no isolated points (it is perfect).

Given a vertex $v\in \V(\cL)$, denote by $B_v$ the corresponding
closed ball. Fix a vertex $v$, and consider the set of balls $B_w$
corresponding to the vertices $w\in S_v$, where
$\Star_{\T(\cL)}(v)=\{[v,w]\ : \ w\in S_v\}$. Then, either $B_w$ are
all disjoint or there exists a $w_0$ such that $B_w\subset B_{w_0}$
for all $w$ and the rest of $B_w$ are disjoint. In fact, if $B_w$
and $B_{w'}$ are not disjoint and both are contained in $B_v$, then
one is inside the other, and hence one is inside the path from $B_v$
to the other, so it does not form and edge. If $B_w$ and $B_{w'}$
contain both $B_v$, then one is contained in the other, and the same
argument applies. Note that the first case happens exactly when
$B_v$ contains $\cL$.

Now, in the first case, we consider the set $\bigcup_{w\in S_v}
B_w$, while in the second case we take
$$\bigcup_{w\in S_v, \ w\ne w_0} B_w \cup \Bc_{w_0}.$$
We are going to see that they are coverings of $\cL$.

Given any point $p\in \cL\cap B_v$, consider a ball
$B(p,\delta)\varsubsetneq B_v$. This ball contains infinite points
of $\cL$ (since $\cL$ has no isolated points), hence we can take
three of them, $p_1=p$, $p_2$ and $p_3\in \cL\cap B(p,\delta)$. The
corresponding vertex $v'=t(p,p_2,p_3)$ is in $\V(\cL)$ and hence
$[v,v']$ contains a vertex $w\in S_v$, and then $p\in B_w$. Now, if
$p\in \cL\setminus B_v$, repeat the same argument with the
complement of an open ball $\Bc(p,\delta)$ such that
$\Bc(p,\delta)\cap B_v=\emptyset$.
\end{Proof}

Finally, we will show that, if $\cL$ is compact and perfect, then
$\cL$ can be identified as the set of ends of $\T(\cL)$. In fact the
bijection is an homeomorphism when considering the ends with the
natural topology.

Recall that a \emph{ray} on a tree $\T=(\V,\E)$ is an infinite
sequence $v_0,v_1,\dots$ of vertices such that $[v_i,v_{i+1}]$ is an
edge and $v_i\neq v_j$ $\forall i\neq j$. Given a progression
$v_0,v_1,\dots,v_n,\dots$ of distinct vertices we say they generate
a ray if the progression formed by the ordered set $\bigcup_{i\ge 0}
\V\cap [v_i,v_{i+1}]$ is a ray, which we call the ray generated by
the $v_n$'s. We denote $\mbox{Rays}(\T)$ the set of rays of $\T$.
Now the ends of a tree $\T$ is the set of equivalence classes of
rays with respect the equivalence relation $\sim$, where
$r=(v_n)\sim s=(w_n)$ if and only if $r\cap s$ is a ray.

The set of ends $\Ends(\T):=\mbox{Rays}(\T)/\sim$ has a natural
topology which has as a subbasis the following sets: for any
oriented (i.e. ordered) edge $e=[v_0,v_1]$, we denote
$$B(e):=\{ r\in \mbox{Rays}(\T) \ | \ e\subset r\}/\sim\subset \Ends(\T)$$
where $e\subset r$ is as ordered sets. Note that $B(\bar e)=B(e)^c$
if $\bar e=[v_1,v_0]$ denotes the opposed edge.

\begin{lemma}
Consider $[v_0,v_1]$ and $[v_1,v_2]$ edges of $\T(L)$, with $v_0\neq
v_2$. Each vertex corresponds to a ball denoted $B_0$, $B_1$ and
$B_2$ respectively. Then if $B_0\cap B_2=\emptyset$ we have that
$B_1$ contains $B_0$ and $B_2$; and if $B_0\cap B_2\neq\emptyset$
either $B_2\subset B_1\subset B_0$ or $B_0\subset B_1\subset B_2$.
\end{lemma}

\begin{Proof}
If $B_0$ and $B_2$ are disjoint but they are linked to the same
vertex then they must be in the ball corresponding to this vertex.
If $B_0\cap B_2\neq\emptyset$, let $p$ be in the intersection. Since
they are linked to $B_1$, $p$ is also in $B_1$, so the possibilities
are either $B_2\subset B_1\subset B_0$ or $B_0\subset B_1\subset
B_2$.
\end{Proof}

A direct consequence of this fact is that if $(B_0,B_1,\dots)$ is a
ray then either there exists $m\geq 0$  such that $B_i\subset
B_{i+1}$ for all $i\geq m$, or $B_{i+1}\subset B_{i}$ for all $i\geq
0$. In fact, once you find two balls in the sequence such that
$B_{j+1}\subset B_j$ then $B_{i+1}\subset B_i$ for all $i\geq j$.

\begin{prop}\label{interseccio}
If, given $m\geq 0$, $B_{i+1}\subset B_i$ for all $i\geq m$, where
the balls corresponds to the vertex of a ray defined as above, then
\begin{equation}\label{aillats1}
    \bigcap_{i\geq m}(B_i\cap\cL)=\{p\} \text{ , where }p\in\cL.
\end{equation}
\end{prop}

\begin{Proof}
For any $i$, $B_i\cap\cL\neq\emptyset$, where $B_i=t(p,p',p'')$ for
$p,p',p''\in\cL$  and two of them are in $B_i$. In fact
$(B_i\cap\cL)\backslash(B_{i+1}\cap\cL)\neq\emptyset$, so since
$B_i\cap\cL$ are closed in $\cL$ and non empty
\begin{equation*}
    \bigcap_{i\geq m}B_i\cap\cL\neq\emptyset.
\end{equation*}

Now we have to see that in the intersection there is a unique point.
Suppose $p_1,p_2$ different in $\bigcap_{i\geq m}B_i\cap\cL$. Then
we take $p_3\in B_m\cap\cL$ different from $p_1$ and $p_2$.
Therefore $v=t(p_1,p_2,p_3)\in \V(\cL)$. But note that
$v_i\in[v_m,v]$ for any $i\geq m$ because $B_m\supset B_i\supset
B(p_1,|p_1,p_2|)=v$. But this is not possible because we know that
$\#[v_m,v]\cap \V(\cL)<\infty$. So $\bigcap_{i\geq
m}B_i\cap\cL=\{p\}$.
\end{Proof}

Conversely, it is clear that a sequence of nested balls $B_i\supset
B_{i+1}$ for all $i>0$ generates a ray when they intersection
$\bigcap (B_i\cap \cL)$ is a point.

\begin{teo}
Let $\cL$ be a compact subset of $P^1(K)$. Then there is a well
defined map
$$\Psi: \mbox{Rays}(\T(\cL)) \to \cL$$
whose image is the set of non isolated points $\cL'\subset \cL$. The
map $\Psi$ determines an homeomorphism between the space of ends and
$\cL'$.
\end{teo}
\begin{Proof}
If $r=(v_1,v_2,\dots)$ with corresponding balls verifying
$B_{i+1}\subset B_i$ for some $i$, then   $B_{i+1}\subset B_i$ for
all $i>m$ and Proposition \ref{interseccio} implies that
$\bigcap_{i\geq m} (B_i\cap \cL)=\{p\}$. We define then $\Psi(r)=p$.
In the other case we have $B_{i+1}\supset B_i$ for all $i$ and we
define $\Psi(r)=\infty$.

If $p\neq\infty$ is in $\mbox{Im}(\Psi)$, so  $p=\Psi(r)$ with
$r=(v_1,v_2,\dots)$, then $v_i=B_i=B(p_{v_i},\rho_{v_i})$ for some
$p_{v_i}\neq p$, $|p-p_{v_i}|=\rho_{v_i}$. Since $B_j\subset B_i$,
for all $j>i>m$, for some $m$, we have $p_{v_i}\neq p_{v_j}$ and
$|p-p_{v_i}|=\rho_{v_i}\rightarrow 0$ when $i$ tends to $\infty$ by
Proposition \ref{interseccio}.

Moreover, any non-isolated point $x\in L$ is in the image of $\Psi$,
since, if $x_i\in L$, $x_i\ne x_j$ for $i\ne j$ and $\lim x_i=x$,
then $v_i:=t(x_1,x_i,x)$ for $i>1$ large enough generate a ray. To
show this, suppose $x\ne \infty$ (the case $x=\infty$ is done by an
analogous argument). Then for $i$ large enough, $v_i$ corresponds to
a ball $B_i$ around $x$ and $B_i\supset B_{i+1}$ since $x_i$
converge to $x$.

Now, it is clear that two rays $r_1$ and $r_2$ have the same image
if they are equivalent, since $\Psi$ only depends of a tail of the
ray. Moreover, if two rays have image $p$, this means that $p$ is
inside the balls of both rays (for large enough index), so they must
be equivalent. That the map $\Psi$ determines an homeomorphism is
clear from the given description.
\end{Proof}

\begin{cor}\label{corvalence}
If $\cL$ is compact and perfect, then $\T(\cL)$ is a locally finite
tree with all vertices of valence strictly bigger than $2$ and
$\Psi$ is surjective.
\end{cor}
\begin{Proof} Only the assertion on the valence of any vertex needs a
comment. Let $v$ be a vertex of $\T(\cL)$, corresponding to a ball
$B(p,\Tog)=t(p,p',p'')$ for some $p,p',p''\in \cL$ with
$\delta=|p-p'|$. Take $\epsilon<\delta$ in $\Tog$. Then, since $\cL$
is perfect, $B(p,\epsilon)\cap \cL$ contains another point $r\in
\cL$, and $t(p,r,p'')\ne t(p,p',p'')$. Similarly, there exists
$p'\ne r'\in \cL\cap B(p',\epsilon)$ and moreover
$B(p',\epsilon)\cap B(p,\epsilon)=\emptyset$, and similarly for
$p''$ (with some minor changes in the case that $p''=\infty$). So we
have vertices $v'$, $v''$ and $v'''$ connected with disjoint paths
to $v$, which means $v$ has valence $3$ or larger.
\end{Proof}

\section{Hyperbolic matrices}\label{sect3.1}

From now on, and for the rest of the paper, we will suppose that $K$
is a field complete with respect to a non-trivial microbial
valuation (if the valuation is non-microbial, there are no
hyperbolic matrices over $K$).

Given any matrix $A\in \GL_2(K)$, we denote by
$\ta(A):=\frac{\tr(A)^2}{\det(A)}\in K$, where $\tr(A)$ is its
trace. It is easily shown that $\ta(A)$ does not depend of the class
in $\PGL_2(K)$, so it gives a well defined map $\ta:\PGL_2(K)\to K$.
Using the natural isomorphism $\AK\cong \PGL_2(K)$, we will use also
$\ta(\gamma)$ for a given $\gamma\in \AK$.

\begin{defi} Given $\gamma\in \AK$, we say that $\gamma$ is \emph{hyperbolic}
if $\ta(\gamma)\in K^*$ and $\ta(\gamma)^{-1}$ is topologically
nilpotent, so a microbe.
\end{defi}

Given any $q\in K^*$, we denote by $\mq\in \AK$ the automorphism
given by $\mq(x)=qx$ for all $x\in K$. The following lemma is a
version of \cite{mumford1972analytic}, lemma 1.1.

\begin{prop}\label{ppp}
Let $\gamma\in \AK$ be any automorphism. Then $\gamma$ is hyperbolic
if and only if there exists $\tau\in \AK$ such that
$\tau\gamma\tau^{-1}=\mq$, where $q\in K$ is a microbe, uniquely
determined by this condition.
\end{prop}
\begin{Proof}
Suppose first that $\gamma=\mq$ with $q$ a microbe. Let $A$ be the
corresponding matrix
$$A=\begin{pmatrix} q&0\\0&1\end{pmatrix}.$$
Then $\tr(A)=q+1$, $\det(A)=q$, $\ta(A)=\frac{(q+1)^2}{q}$, so
$|\ta(A)|=|q|^{-1}$ since $|q+1|=1$. It is clear that $\ta(\gamma)$
is invariant by conjugation (as they are the determinant and the
trace), which shows one implication.

Now suppose $\gamma$ is hyperbolic, with associated matrix $A\in
\GL_2(K)$. We take a representative $A$ with coefficients in
$\mathcal{O}_K$. Let $f(x)=x^2-ax+b$ be the characteristic
polynomial of $A$, so $\ta(A)=\frac {a^2}{b}$. By hypothesis
$t=\frac{|f(0)|}{|f'(0)|^2}$ is a microbe. Using Hensel's Lemma
(which holds since $K$ is complete) we see that there exists
$\alpha\in\mathcal{O}_K$ with $f(\alpha)=0$. Therefore there exist
also $\beta\in\mathcal{O}_K$ such that $f(x)=(x-\alpha)(x-\beta)$.
Note that $\beta\in\mathcal{O}_K$ because $\alpha$ and
$\alpha+\beta\in\mathcal{O}_K$.

Moreover $|\alpha|\ne |\beta|$, since by hypothesis $$
\left|\frac{\alpha\beta}{(\alpha+\beta)^2}\right|<1 ,$$ so
$|\alpha|=|\beta|$ implies $|\alpha|^2<|2\alpha|^2\le |\alpha|^2$
which is a contradiction. Summarizing we have that or
$|\frac{\alpha}{\beta}|<1$ or $|\frac{\beta}{\alpha}|<1$, hence one
of these is equal to $|\ta(A)^{-1}|$, therefore a microbe.

The unicity assertion follows easily since $\mu_q$ and $\mu_{q'}$
are conjugate if and only if $q=q'$.
\end{Proof}

Given an hyperbolic automorphism $\gamma\ne \id$, we denote by
$q_{\gamma}\in K^*$ the unique microbe such that $\gamma$ is equal
to $\mq$ modulo conjugation. We denote also
$\varrho(\gamma)=|q_{\gamma}|\in \Tog$.

\begin{cor}\label{cor4.10}
Let $\gamma\in \AK$ be an hyperbolic automorphism. Then $\gamma$ has
two fixed points defined over $K$. For any $p\in \bP^1(K)$ such that
$\gamma(p)\ne p$, the limits $\lim_{n\to \infty}\gamma^n(p)$ and
$\lim_{n\to -\infty}\gamma^n(p)$ exists and are equal to the two
fixed points by $\gamma$.\end{cor}

\begin{Proof}
By conjugation, we reduce to the case $\gamma=\mq$ with $q$
topologically nilpotent. In this case $\lim_{n\to
\infty}|q^n||p|=|0|$ so $\lim_{n\to \infty}\gamma^np=0$ and the
fixed points are $0$ and $\infty$.
\end{Proof}

Hyperbolic automorphisms are specially interesting for us since they
don't fix any element of $\BK$.

\begin{lemma}\label{lemma1}
 Let $\gamma\ne \id$ be an hyperbolic automorphism. Then $t(p_1,p_2,p_3)\neq t(\gamma p_1,\gamma p_2,\gamma p_3)$ for all $(p_1,p_2,p_3)\in \wK$.
\end{lemma}
\begin{Proof}
By conjugation, we can and will suppose that $\gamma=\mq$, for $q$
topologically nilpotent. So $t(p_1,p_2,p_3)$ corresponds (reordering
if is necessary) to $B(p_1,|p_1-p_2|)$ but also $t(qp_1,qp_2,qp_3)$
corresponds to $B(qp_1,|q||p_1-p_2|)$ and since $|q|<1$ one has
$|q||p_1-p_2|<|p_1-p_2|$ so $B(p_1,|p_1-p_2|)\supsetneq
B(qp_1,|q||p_1-p_2|)$.
\end{Proof}

In the following we will show some properties that characterize
hyperbolic matrices. Recall that a $\gamma \in \Aut(K)$, with
$\gamma\ne \id$, either has one or two fixed points over the
algebraic closure $\bar K$. The hyperbolic automorphisms have two
fixed points, defined over $K$.

The following result translates the definition of hyperbolic
automorphisms to the action on the $\Lambda$-tree $\BK$.

\begin{prop}\label{hyperbolicballs} If for a given $\vp\in\AK$ there
exist a ball $B$ such that $\vp(B)\subset B$ or $B\subset \vp(B)$
and it verifies that $d(B,\vpB{\vp}{B})^{-1}\in \Tog$ is
topologically nilpotent, then the automorphism $\vp$ is hyperbolic
and $\varrho(\vp)=d(B,\vpB{\vp}{B})^{-1}$.

Moreover, for any hyperbolic automorphism $\vp$ there exist such a
ball $B$ for $\vp$. Even more, for any ball $B$,
$d(B,\vpB{\vp}{B})^{-1}$ is topologically nilpotent.
\end{prop}

\begin{Proof} We denote by $\rho:=d(\vpB{\vp}{B},B)$, by $\delta=\varrho(B)$
and we fix $p\in B$, so $B=B(p,\delta)$. First of all we are going
to proof that $\vp$ has two fixed points, $p_0\in K$ and
$p_{\infty}\in \bP^1(K)$. We have three cases to consider: either
$\vp(B)\subset B$, or $\vp(B)$ is a ball and $B\subset \vp(B)$, or
$\vp(B)$ is the complement of a ball.

First we show that in the first case $B$ contains a point $p_0$
fixed by $\vp$. By applying Lemma \ref{pairofballs} to
$\vp(B)\subset B$ and $\vp$ we get by induction that
$\rho=d(\vp^{n+1}(B),\vp^n(B))$ for all $n\ge 1$. Hence for all
$m\ge 1$ we have $\vp^m(B)=B(\vp^m(p),\rho^m\delta)$, and
$\rho^m\delta\to 0$ for $m\to \infty$, so $\bigcap_{n\ge
1}\vp^n(B)=\{p_0\}$ where $p_0$ is a point necessarily fixed by
$\vp$ (in fact $p_0=\lim_{n\to \infty} \vp^n(p)$). Note that then
$\rho=|\vp'(p_0)|$.

Now consider the sets $\vp^{-n}(B)$ for $n\ge 1$. We have two
options: either $\vp^{-n}(B)$ are balls for all $n\ge 1$, or there
exists one $n_0\ge 1$ such that $\vp^{-n_0}(B)$ is the complement of
a ball. If the first option happens, then we necessarily have that
$\vp^n(B)\supset \vp^{n+1}(B)$ for all $n\in \bZ$ and
$\varrho(\vp^n(B))=\rho^n\delta$. Hence $\bigcup_{n\in \bZ}
\vp^n(B)=K$ and then $\infty$ must be a fixed point of $\vp$, which
shows the result in the first case under this hypothesis.

If it happens the second option, we can and will suppose that
$n_0=1$, and then $\vp^{-1}$ and $B$ are in the third case above, so
$\vp^{-1}(B)$ is the complement of a ball, a case we will consider
now.

So, we suppose that we are in the third case, thus $\infty\in
\vp(B)\supset B$. Therefore
$\vpB{\vp}{B}=B(\vp(\infty),|\vp'(\infty)|\delta^{-1})$ by
Proposition \ref{Actionont} (2).

First observe that $\vpB{\vp}{B}\cap B=\emptyset$. To show this note
that, if $\vpB{\vp}{B}\cap B \ne \emptyset$, then either
$\vpB{\vp}{B}\varsubsetneq B$ or $B\varsubsetneq \vpB{\vp}{B}$,
since they cannot be equal as they are at distance larger than one.
But in the first case some points of $B$ will not be in $\vp(B)$,
while in the second $\vp(B)\cap B=\emptyset$.

But then we have that $\vp(\vpB{\vp}{B})\subset \vpB{\vp}{B}$ and
that $\vp^{-1}(B)\subset B$. This is because an element of
$\vp(\vpB{\vp}{B})$ which is not in $\vpB{\vp}{B}$ will be then on
$\vp(B)$ (since $\vpB{\vp}{B}\cup \vp(B)=\bP^1(K)$), hence will be
the image by $\vp$ of an element in $B\cap \vpB{\vp}{B}=\emptyset$.
And the same argument shows that an element of $\vp^{-1}(B)$ not in
$B$ is in $\vp^{-1}(\vpB{\vp}{B})$, which again cannot happen.

Hence, by applying the first result under the fist case, we have
that $\vp$ has a fixed point $p_0\in \vpB{\vp}{B}$, and $\vp^{-1}$,
and hence $\vp$, has a fixed point $p_{\infty}$ in $B$, which shows
the result in the third case.

Therefore, and returning now to the first case, and when $n_0=1$, we
have that $\vp^{-1}$ has a fixed point $p_{\infty}$ in
$\vpB{\vp^{-1}}{B}$, which also finish with the first case.

Finally, in the remaining second case, so $\vp(B)\supset B$ and it
is a ball, then $\vp^{-1}$ is in the first case, so it has two fixed
points, and thus also $\vp$.

Take now $\tau\in \AK$ such that $\tau(p_0)=0$ and
$\tau(p_{\infty})=\infty$, i.e. $\tau(t)=(t-p_0)/(t-p_{\infty})$ (or
$\tau(t)=t-p_0$ if $p_{\infty}=\infty$). Then $\mu=\tau\circ\vp\circ
\tau^{-1}$ has fixed points $0$ and $\infty$, hence it is of the
form $\mu=\mu_q$ for some $q\in K$. One shows easily that
$|q|=\rho$, which proofs that $\vp$ is hyperbolic by proposition
\ref{ppp} and that $\varrho(\vp)=\rho$.

The assertion about the existence of such a ball for a given
hyperbolic $\vp$ is easy. Let $\tau\in \AK$ be such that
$\tau\gamma\tau^{-1}=\mq$, where $q\in K$ is topologically
nilpotent. Suppose that $\tau^{-1}(0)\ne \infty$, and consider a
ball $B_0$ with center $0$ such that $\tau(\infty)\notin B_0$. Then
$B=\tau^{-1}(B_0)$ is a ball such that $\vp(B)\subset B$ and
verifies the hypothesis. And in case $\tau^{-1}(0)=\infty$, one does
the construction for $\vp^{-1}$, getting at the end $B\subset
\vp(B)$.

The final assertion can be shown just in the case $\vp=\mq$, and in
this case it is easy. In this case, if $0\in B$, then
$d(B,\vpB{\vp}{B})=\varrho(\vp)$, and if $0\notin B=B(p,\delta)$,
then
$$d(B,\vpB{\vp}{B})=\varrho(\vp)d(B,B(p,|p-0|))d(\vp(B),B(\vp(p),|\vp(p)-0|)).$$
\end{Proof}

In fact, it is easy to construct hyperbolic automorphisms verifying
the hypothesis of the lemma for a given balls.

\begin{lem}\label{twoballs} Given any balls $B'$ and $B$ with
$d(B',B)^{-1}$ topologically nilpotent, there exists $\vp\in\AK$
hyperbolic with $\varrho(\vp)=d(B',B)^{-1}$ and such that
$\vpB{\vp}{B}=B'$.
\end{lem}

\begin{Proof} Choose $q\in K^*$ such that $|q|=d(B',B)^{-1}$.
After reordering the balls, we can suppose that, either $B\cap
B'=\emptyset$ or $B'\subset B$. Consider an automorphism $\tau$
sending a point of $B'$ to 0, and, in case $B\cap B'=\emptyset$, a
point of $B$ to $\infty$. Then the hyperbolic automorphism we are
looking for is $\vp:=\tau^{-1}\mq\tau$.
\end{Proof}

The following result shows that that the hiperbolicity of $\gamma$
is directly related to the compacity of the closure of the orbit of
any point by $\gamma$.

\begin{teo}\label{diagon}
Suppose $\gamma\in \Aut(K)$ has two fixed points over the algebraic
closure $\bar K$. Let $\Gamma=\langle \gamma\rangle$ be the subgroup
generated by $\gamma$. For any $p\in\bP^1(K)$, let
$\Gamma_p=\{\gamma^np \ | \ n\in\mathbb{Z}\}$ be the orbit of $p$.
Then
\begin{equation*}
    \overline{\langle\gamma\rangle p}\text{ is compact for all }p\in\bP^1(K)
    \Leftrightarrow\begin{cases}\gamma \text{ is hyperbolic, } \\ \gamma\text{ is of finite
    order.}\end{cases}
\end{equation*}
\end{teo}

\begin{Proof}
The reverse implication is easy. First note that if $\gamma$ is of
finite order it means that $\langle\gamma\rangle p$ is finite, so it
is compact. It remains to show that if $\gamma$ is hyperbolic then
$\overline{\langle\gamma\rangle p}$ is compact. Since an hyperbolic
matrix $\gamma$ is conjugated to $\mq$ given by $\mq(x)=qx$ with $q$
topologically nilpotent, the closure of the orbit of $\gamma$ will
have the same topological structure as $\mq$. Note that
\begin{equation*}
   \langle\mq\rangle p=\begin{cases}\{p\} \text{ if }p=0 \text{ or }p=\infty, \\ \{q^np \ | \ n\in\mathbb{Z}\text{ if }p\neq0,\infty\}.\end{cases}
\end{equation*}
We have to see that the second case is also compact. First note that
$\overline{\langle\gamma\rangle p}=\{q^np \ | \
n\in\mathbb{Z}\}\cup\{0,\infty\}$ because
$\lim_{n\rightarrow\infty}q^np=0$ and
$\lim_{n\rightarrow-\infty}q^np=\infty$. Let $
\overline{\langle\gamma\rangle p}=\bigcup_{i\in I}B_i$ be a covering
of open balls. Since $0$ is in the covering we have a ball that
contains it, which must be of the form $B=B(0,\rho)$. By the same
reasoning there is a ball that contains $\infty$ like
$B'=\Bc(0,\rho')$. But then $\overline{\langle\gamma\rangle
p}\backslash(B\cup B')$ is finite, since there exists $n_0,n_0'\ge
1$ such that $q^np\in B(0,\rho)$ for all $n\geq n_0$ and $q^{-m}p\in
\Bc(0,\rho')$ for all $m\geq n'_0$.

To see the direct implication, observe first that, by extending the
field to an at most a degree two extension, we can reduce to the
case that $\gamma$ has two fixed points defined over $K$ (as being
hyperbolic or of finite order is independent of the field). Hence we
are reduced to show that the automorphism $\mq$, with $|q|\le1$ and
$q$ not topologically nilpotent, does have a non compact orbit. We
will consider the orbit of $1$, i.e. $\Gamma 1=\{q^n\ | \ n\in
\bZ\}$,and show that $\overline{\Gamma 1}=\Gamma 1$ and it is not
compact. Since $q$ is not topologically nilpotent, there exists
$\lambda\in \Tog$ such that $\lambda^{-1}>|q^n|>\lambda>0$ for all
$n\in \bZ$. So, for $n>m$,
\begin{equation*}
    |q^n-q^m|=|q|^m|q^{n-m}-1|=|q|^m|1|>\lambda.
\end{equation*}
So $\Gamma1$ is covered by $ \bigcup_{n\in\mathbb{Z}}B(q^n,\lambda)$
which satisfies
$$
    B(q^n,\lambda)\cap \overline{\{\gamma^n(1) \ | \ n\in\mathbb{Z}\}}=\{q^n\}
$$ and are pairwise disjoint, so we can not remove any
ball.  Moreover $\Gamma1$ is closed, since it has no accumulation
point in $K$, and, because $|q|^{-n}<\lambda^{-1}$ for all
$n\in\mathbb{Z}$, so we have $\Gamma1\subset B(0,\lambda^{-1})$,
$\infty$ is not an accumulation point.
\end{Proof}

We can remove the condition of having two fixed points if we exclude
the so called mixed characteristic case.

\begin{teo}
Let $K$ be a complete field and algebraically closed and either
$\mbox{char}(K)=p>0$ or
$\mbox{char}(\mathcal{O}_K)/\mathfrak{m}_K=0$.  Then for $\gamma\in
PGL_2(K)$
\begin{equation*}
    \overline{\langle\gamma\rangle p}\text{ is compact for all }p\in\bP^1(K)\Leftrightarrow\begin{cases}\gamma \text{ is hyperbolic, }
    \\ \gamma\text{ is of finite order.}\end{cases}
\end{equation*}
\end{teo}
\begin{Proof}
Since $K$ is algebraically closed we can suppose that either
$\gamma$ is diagonalizable, a case already done in theorem
\ref{diagon}, or it is conjugate to $\psi_a(t)=t+a$ for some $a\in
K\backslash\{0\}$. We can and will suppose $\phi=\psi_a$.

Consider then the orbit of $0$, $\langle\gamma\rangle 0=\{na \ | \
n\in\mathbb{Z}\}$. If $\mbox{char}(K)=p>0$, then $\gamma$ has finite
order equal to $p$. Now, if $\mbox{char}(K)=0$ and
$\mbox{char}(\mathcal{O}_K/\mathfrak{m}_K)=0$, then $|n|=1$ for all
$n\in\bZ$ so $|na|=|a|$ for all $n\in\mathbb{Z}$. Moreover if $n\neq
m$ then $|na-ma|=|n-m||a|=|a|$. Then all points in the bounded set
$\langle\gamma\rangle 0$ are isolated, so the set is closed, but it
has an infinite number of points, so it is not compact.
\end{Proof}

\begin{remark}
The result is false if $\mbox{char}(K)=0$ and
$\mbox{char}(\mathcal{O}_K/\mathfrak{m}_K)=p>0$, because in this
case $|.|\big|_{\mathbb{Q}}=|.|_p^\epsilon$, where $|.|_p$ denotes
the $p$-adic absolute value and $\epsilon\in \bR_{>0}$. Then the
closure of the orbit of $0$ for $\psi(t)=t+1$ is $\mathbb{Z}_p$, the
$p$-adic integers, which is compact. For any $p\in\bP^1(K)$ the
closure of the orbit is a translate of $\bZ_p$, hence compact as
well.
\end{remark}

\section{Schottky Groups}\label{sect3}

Now we are going to define an analogous of Schottky groups in our
context, as a finitely generated subgroup, not cyclic, formed by
hyperbolic matrices, plus a compacity condition.

\begin{defi}
For $\Gamma\subset \AK$ a subgroup, we define
\begin{equation*}
    \Fix(\Gamma)=\{p\in\bP^1(K) \ | \ \exists\gamma\in\Gamma,\gamma\neq \mbox{Id} \text{ with } \gamma(p)=p\}
\end{equation*}
and $\cL_\Gamma=\overline{\Fix(\Gamma)}$ its closure.
\end{defi}

\begin{defi}
Let $\Gamma\subseteq PGL_2(K)$ be a subgroup. We say that it is a
\emph{Schottky group} if
\begin{enumerate}
    \item $\Gamma$ is finitely generated
    \item every element of $\Gamma$ different from the identity is hyperbolic
    \item $\overline{\Gamma P}$ (the closure of the orbit of $P$) is compact for all $P\in\bP^1(K)$.
    \item $\Gamma$ is not cyclic.
\end{enumerate}
\end{defi}

Note that a Schottky group is torsion free. Note also that a
finitely generated but not cyclic subgroup of a Schottky group is a
Schottky group.

We will show that for any Schottky group $\Gamma$, the set
$\cL_{\Gamma}$ is compact and perfect.

\begin{lemma}\label{sublemma1}
If $q$ and $r\in K^*$ topologically nilpotent,
the subgroup $\Gamma$ that they generate does not have torsion
elements and $q^m=r^n$ for some $n$ and $m\in
\mathbb{Z}\setminus\{0\}$, then $\Gamma$ is cyclic.
\end{lemma}
\begin{Proof}
We can suppose $m$ and $n$ are coprime, since, if $s=\gcd(n,m)$,
$n=sn'$, $m=sm'$, then $(q^{n'}r^{-m'})^s=q^{n}r^{-m}=1$ so
$q^{n'}r^{-m'}=1$ since $\Gamma$ has no elements of finite order.
But now, there exists $a,b\in \mathbb{Z}$ with $am+bn=1$ and we have
$(q^{b}r^{a})^m=q$ and $(q^{b}r^{a})^n=r$, so $q$ and $r$ belong to
the subgroup generated by $q^br^a$.
\end{Proof}

\begin{lemma}\label{sublemma2} If $q$ and $r\in K^*$ are topologically nilpotent,
the subgroup $\Gamma$ that they generate does not have torsion
elements and $|q|^m\ne |r|^n$ for all $n$ and $m\in
\mathbb{Z}\setminus\{0\}$, then $\mathcal{W}:=\overline{\{q^nr^m \ :
\ (n,m)\in \mathbb{Z}^2\}}$ is not compact.
\end{lemma}
\begin{Proof} We will suppose that $1>|r|>|q|$. We will show that
$\mathcal{W}$ contains infinitely many isolated points contained in
a bounded set. Consider
$$W:=\mathcal{W}\cap \{x\in K \ : \ 1\ge |x|>|q|\}$$
Now, for any $x\in W$, take the ball $B(x,|q|)$. Observe that for
any $y\in B(x,|q|)$, $|y|=|(y-x)+x|=\max(|y-x|,|x|)=|x|$, since
$|y-x|\le q<|x|$. But by hypothesis no two elements in $\mathcal{W}$
have the same valuation, hence $\mathcal{W}\cap B(x,|q|)=\{x\}$ for
any $x\in W$.

But the set $W$ contains an infinite number of points, since, that
for any $m\ge 1$, there exists $f(m)\in \bZ$ such that $|r|^{f(m)}>
|q|^m \ge |r|^{f(m)+1}$, hence $x_m:=q^mr^{-f(m)}$ is in $W$ for all
$m$, since $1>|q|^m|r|^{-f(m)}\ge |r| >|q|$.

Then
$$\mathcal{W}=\overline{\{q^nr^m \ : \ (n,m)\in\mathbb{Z}^2\}}\subset \left(\bigcup_{w\in W} B(w,|q|) \right) \cup B(0,|q|) \cup B^0(0,1)^c $$
and no ball can be removed.
\end{Proof}

\begin{lemma}\label{fixof2}
Let $\Gamma$ be a Schottky group. Then, for any $id\ne \gamma$ and
$\tau \in\Gamma$, either $\Fix(\gamma)=\Fix(\tau)$, and then
$\gamma$ and $\tau$ belong to a cyclic subgroup, or
$\Fix(\gamma)\cap \Fix(\tau)=\emptyset$ .
\end{lemma}

\begin{Proof} If $\Fix(\gamma)=\Fix(\tau)=\{p_0,p_1\}$, then, after
conjugation, we can suppose $\gamma(x)=qx$ and $\tau(x)=rx$ for some
$q$ and $r\in K^*$ topologically nilpotent. We will see that then
$q^m=r^n$ for some $n,m\in \bZ\setminus\{0\}$, and then lemma
\ref{sublemma1} implies that $\gamma$ and $\tau$ belong to the same
cyclic subgroup. If $q^m\ne r^n$ for all $n$ and $m\in
\mathbb{Z}\setminus\{0\}$, but $|q|^m=|r|^n$ for some  $n$ and $m\in
\mathbb{Z}\setminus\{0\}$, then $q^mr^{-n}$ is not $1$ but has
valuation $1$. Hence $\gamma^m\tau^{-n}\in \Gamma$ and it is not
hyperbolic. And if  $|q|^m\ne |r|^n$ for all $n$ and $m\in
\mathbb{Z}\setminus\{0\}$, Lemma \ref{sublemma2} gives us a
contradiction.

If $\Fix(\gamma)\cap \Fix(\tau)=\{p_0\}$ for some point $p_0$, we
can reduce to the case, after conjugation, that $p_0=\infty$,
$\Fix(\gamma)=\{0,\infty\}$ and $\Fix(\tau)=\{1,\infty\}$. Then
 $\gamma(x)=qx$ and $\tau(x)=rx+(1-r)$  for some
$q$ and $r\in K^*$ topologically nilpotent. But then
$\gamma\tau\gamma^{-1}\tau^{-1}(x)=x-(1-q)(1-r)$, which is clearly
non hyperbolic.
\end{Proof}

\begin{lemma}\label{orbitoffixed}
Suppose $\Gamma$ is a Schottky group. Consider $p\in \cL_{\Gamma}$.
Then $ \overline{\Gamma p}=\cL_{\Gamma}$.
\end{lemma}
\begin{Proof}
If $p$ is fixed by $\gamma\ne 1\in \Gamma$, then $\alpha (p)$ is
fixed by $\alpha^{-1}\gamma\alpha$ for any $\alpha\in \Gamma$. So
$\overline{\Gamma p}\subset\cL_{\Gamma}$.

Now, if $p'$ is another point in $\cL_{\Gamma}$, with $\gamma(p')\ne
p'$, and fixed by some $\alpha\in \Gamma$, by the previous lemma
\ref{fixof2}, then $\alpha(p)\ne p$ and hence $\alpha^n(p)\to p'$
for $n\to \pm \infty$ by corollary \ref{cor4.10}. So $p'\in
\overline{\Gamma p}$. So all points fixed by some $\alpha\in\Gamma$,
except may be the other point different from $p$ fixed by $\gamma$,
are in $\overline{\Gamma p}$, which imply that its closure, which is
$\cL_{\Gamma}$, is contained in $\overline{\Gamma p}$.

Finally, if $p\in L_{\Gamma}$ is the limit of points $p_n$ fixed by
some $\gamma_{n}\in \Gamma$, then any point in $\overline{\Gamma p}$
is limit of points in $\overline{\Gamma p_n}=\cL_{\Gamma}$, so it is
in $\cL_{\Gamma}$. The reverse inclusion is also clear.
\end{Proof}

\begin{teo}\label{compactandperfect}
Suppose $\Gamma$ is a Schottky group. Then the set $\cL_{\Gamma}$ is
perfect and compact.
\end{teo}
\begin{Proof}
It is compact since, by the previous lemma \ref{orbitoffixed},
$\cL_\Gamma=\overline{\Gamma p}$ for some $p\in\cL_\Gamma$, and
$\overline{\Gamma p}$ is compact by definition of Schottky group.

Let $p\in \bP^1(K)$ be fixed by $\gamma\in \Gamma$. Take $p'\in \cL$
not fixed by $\gamma$ (for example, fixed by some $\gamma'$ not
contained in the subgroup generated by $\gamma$, that it exists
because $\Gamma$ is not cyclic). Then $\gamma^n(p') \to p$ when
$n\to \infty$ or when $n\to -\infty$. Hence no point fixed by some
$\gamma\ne 1$ in $\Gamma$ is isolated, so the same is true for the
points in the closure.
\end{Proof}

\section{The finite graph associated to a Schottky group.}

The main aim of this section is to show that the quotient by
$\Gamma$ of the tree associated to $\cL_{\Gamma}$ for a Schottky
group $\Gamma$ is finite, and that the quotient map is the universal
cover, hence identifying $\Gamma$ with the fundamental group. We
will denote $\cT_\Gamma:==\T(\cL_\Gamma)$.

\begin{teo}\label{main} Let $\Gamma$ be a Schottky group on a field complete
with respect to a valuation. Then the tree $\cT_\Gamma$ is locally
finite, the group $\Gamma$ acts freely on $\mathcal{T}_\Gamma$ and
the quotient $G_\Gamma:=\mathcal{T}_\Gamma/\Gamma$ is a finite
graph.
\end{teo}

We will prove the theorem along the section. The first part of the
result is a consequence of theorem \ref{compactandperfect} and the
results of Section \ref{sect2}. The group acts freely because of
Lemma \ref{lemma1}, which says that for all $\gamma\in\Gamma$
different from the identity and for all $v\in
\V(\mathcal{T}_\Gamma)$, $\gamma(v)\neq v$.

So we can take the quotient $G_\Gamma=\mathcal{T}_\Gamma/\Gamma$ and
the quotient map $\mathcal{T}_\Gamma\rightarrow G_\Gamma$ is the
universal cover. We only need to show that the graph $G_\Gamma$ is
finite.

\begin{defi}
Let $B_\Gamma\subset\Gamma$ be a finite set of generators verifying
that, if $\gamma\in B_\Gamma$, then $\gamma^{-1}\in B_\Gamma$, and
$\id\in B_\Gamma$. For a fixed vertex $\omega\in\mathcal{T}_\Gamma$
we consider $S_\omega=\{\gamma\omega \ | \ \gamma\in B_\Gamma\}$,
which is a finite set of vertices. We denote
$\mathcal{T}_{S_\omega}=\bigcup_{v_1,v_2\in
S_\omega}[v_1,v_2]=\bigcup_{\gamma\in
B_\Gamma}[\omega,\gamma\omega]$, the minimal finite subtree that
contains $S_\omega$. Finally we denote
\begin{equation*}
    \cT_{B_\Gamma,\omega}=\bigcup_{\gamma\in\Gamma} \gamma(\mathcal{T}_{S_\omega}).
\end{equation*}
\end{defi}

Our aim will be to show in a series of lemmata that
$\cT_{B_\Gamma,\omega}=\cT_{\Gamma}$, and the finiteness of
$G_{\Gamma}=\mathcal{T}_{B_\Gamma,\omega}/\Gamma$ follows. This is
because $\mathcal{T}_{B_\Gamma,\omega}/\Gamma$ has a finite number
of vertices, since
\begin{equation*}
    \V(\cT_{S_\omega})\rightarrow \V(\cT_{B_\Gamma,\omega}/\Gamma)=\V(\cT_\Gamma/\Gamma)
\end{equation*} and $\V(\mathcal{T}_{S_\omega})$ is finite, and the tree
$\cT_\Gamma$ is locally finite, hence also $G_\Gamma$.

\begin{lemma}\label{lemma2}
\
 \begin{enumerate}
 \item $\forall\gamma\in\Gamma$, $[\omega,\gamma\omega]\subset\mathcal{T}_{B_\Gamma,\omega}$,
 \item $\forall \gamma\neq\gamma'\in\Gamma$, $[\gamma\omega,\gamma'\omega]\subset\mathcal{T}_{B_\Gamma,\omega}$.
 \end{enumerate}
\end{lemma}

\begin{Proof}
Since $\gamma\in\Gamma$, then
$\gamma=\gamma_1\gamma_2\dots\gamma_n$, where $\gamma_i\in
B_{\Gamma}$. Then
\begin{equation*}
[\omega,\gamma\omega]\subset[\omega,\gamma_1\omega]\cup[\gamma_1\omega,\gamma_1\gamma_2\omega]\cup\dots
\cup[\gamma_1\gamma_2\dots\gamma_{n-1}\omega,\gamma_1\gamma_2\dots\gamma_{n}\omega]
\end{equation*}
and also each
\begin{equation*}
    [\gamma_1\dots\gamma_{i}\omega,\gamma_1\dots\gamma_{i+1}\omega]\subset(\gamma_1\dots\gamma_i)(\mathcal{T}_{B_\Gamma,\omega})=\mathcal{T}_{B_\Gamma,\omega}.
\end{equation*}
To show the second part, we can divide the path as follows
$[\gamma\omega,\gamma'\omega]\subset [\gamma\omega,\omega]\cup
[\omega,\gamma'\omega]\subset\mathcal{T}_{B_\Gamma, \omega}$, or we
also can argue that
$[\gamma\omega,\gamma'\omega]=\gamma[\omega,\gamma^{-1}\gamma'\omega]\subset\gamma(\mathcal{T}_{B_\Gamma,
\omega})=\mathcal{T}_{B_\Gamma, \omega}$.
\end{Proof}

\begin{lemma}\label{subtree}
The graph $\mathcal{T}_{B_\Gamma,\omega}$ is connected, hence it is
a subtree.
\end{lemma}

\begin{Proof}
Consider two vertices $v_1$ and
$v_2\in\mathcal{T}_{B_\Gamma,\omega}$. Then $v_1=\gamma_1(\omega_1)$
and $v_2=\gamma_2(\omega_2)$ for some $\omega_1,\omega_2\in
T_{B_\Gamma,\omega}$ and some $\gamma_i\in \Gamma$ for $i=1,2$.
Since $T_{B_\Gamma,\omega}$ is connected, there exist paths
$[\omega_1,\omega]$ and $[\omega_2,\omega]$ in $T_{B_\Gamma,\omega}$
and from this one has that
$\gamma_1[\omega_1,\omega]=[v_1,\gamma_1(\omega)]$ and
$\gamma_2[\omega_2,\omega]=[v_2,\gamma_2(\omega)]$ are contained in
$\mathcal{T}_{B_\Gamma,\omega}$. By Lemma $\ref{lemma2}$ we have
$[\gamma_1(\omega),\gamma_2(\omega)]\subset\mathcal{T}_{B_\Gamma,\omega}$,
so
\begin{equation*}
    [v_1,v_2]\subset [v_1,\gamma_1(\omega)]\cup[\gamma_1(\omega),\gamma_2(\omega)]\cup[\gamma_2(\omega),v_2] \subset\mathcal{T}_{B_\Gamma,\omega}.
\end{equation*}
\end{Proof}


\begin{lemma}\label{minimalsubtree}
Let $\Gamma$ be a Schottky group. Let $\mathcal{T}'\subset
\mathcal{T}_\Gamma$ be a non-empty subtree which is invariant by
$\Gamma$. Then $\mathcal{T}'=\mathcal{T}_\Gamma$.
\end{lemma}
\begin{Proof}
First, $\mathcal{T}'$ is infinite since it contains an infinite
number of vertices: the ones of the form $\gamma(v)$, for some $v\in
\mathcal{T}'$ and $\gamma\in \Gamma$.

Let $\cL'$ be the image of $\mathcal{T}'$ with respect to the map
$$\Psi: \mbox{Rays}(\T(\cL_{\Gamma})) \to \cL_{\Gamma}.$$
Clearly $\cL'$ is invariant by $\Gamma$, and non-empty since $\cT'$
is infinite, so it contains some ray. Take $p\in \cL'$. Then $\Gamma
p \subset \cL'$. By lemma \ref{orbitoffixed} we have
$\cL_{\Gamma}=\overline{\Gamma p}\subset \overline{\cL'}\subset
\cL_{\Gamma}$, thus $\overline{\cL'}=\cL_{\Gamma}$

Now, observe that for any $x$ and $y\in \cL'$, all the points of the
form $t(x,y,z)$, for $z\in \cL$, are in fact in $\mathcal{T}'$. To
show this, observe that $x\in \cL'$ implies that the ray
$[t(x,y,z),x]$ contains some vertex $v_x$ of $\mathcal{T}'$ (in
fact, infinitely many). The same happens for $y$, so $[t(x,y,z),y]$
contains a vertex $v_y$ of $\mathcal{T}'$. But $t(x,y,z)\in
[v_x,v_y]\subset \mathcal{T}'$ since $\mathcal{T}'$ is a tree, hence
connected.

But $\cL'$ is closed. Effectively, suppose we have a progression of
distinct points $p_n\in \cL'$ such that $p_n\to p\in \cL$ when $n\to
\infty$. Then the vertices $v_i:=t(p_1,p_2,p_i)$ for $i>2$ are in
$\mathcal{T}'$, and they generated a ray $r$. Then $\Psi(r)=p$, and
hence $p\in \cL'$. So $\cL_{\Gamma}=\overline{\cL'}=\cL'$, and hence
$\cT'=\cT_\Gamma$.
\end{Proof}

As a consequence, we can finish the proof of the theorem \ref{main}.
We have $\mathcal{T}_{B_\Gamma,\omega}$ is invariant by $\Gamma$ by
definition and it is a subtree by corollary \ref{subtree}, so
$\mathcal{T}_{B_\Gamma,\omega}=\mathcal{T}_\Gamma$ by lemma
\ref{minimalsubtree}.

\section{Explicit constructions and fundamental domains.}

The aim of this section is to construct explicit examples of
Schottky groups as well as to give a criterium to decide if some
subgroup of $\AK$ is Schottky.

Recall from Lemma \ref{twoballs} that given two balls $B$ and $B'$
such that $d(B,B')^{-1}$ is topologically nilpotent, there exists an
hyperbolic automorphism $\vp$ such that $B'= \vpB{\vp}{B}$ and
$\varrho(\vp)=d(B,B')^{-1}$. This automorphism depends on the
election of centers of $B$ and $B'$ and of a $q\in K^*$ such that
$|q|=d(B,B')^{-1}$. If $B$ and $B'$ are disjoint, we can even take
$\vp$ such that $\infty \in\vp(B)$, so it is not a ball.

The proof of the following theorem is essentially the same as in
part (4.1.3) in \cite{GvdP}, with the appropriate modifications. We
tried, however, to write again all details of the proof in order to
clarify it.

\begin{teo}\label{fundamental} Take $g\ge 2$ an integer. Suppose we are given $2g$ pairwise
disjoint balls $B_1,\dots,B_{2g}$ such that
$\rho_{i,j}:=d(B_i,B_j)^{-1}\in \Tog$ are topologically nilpotent
elements for all $i\ne j\in\{1,\dots,2g\}$. Let $\vp_i\in \AK$ be
hyperbolic automorphisms such that $B_{i+g}=\vpB{\vp_i}{B_i}$ and
$\infty\in \vp_i(B_i)$ for all $i\in\{1,\dots,g\}$.

Then $\Gamma:=\langle\vp_1,\dots,\vp_g\rangle \subset \AK$ is a
Schottky subgroup.
\end{teo}

To proof the compactness condition necessary for the theorem we will
use the following lemma, which is a version of the Heine-Borel
property of compact sets for $\Tog$-metric spaces. The proof is
almost the same as the classical result.

\begin{lem}\label{Heine-Borel} Suppose $W\subset K$ is a subset such that, for
every $\delta\in \Tog$ there exists a finite subset
$S_{\delta}\subset \BK$ of balls such that $\varrho(B)\le \delta$
for any $B\in S_{\delta}$ and $W\subset \bigcup_{B\in S_{\delta}}
B$. Then $\overline{W}$ is compact.
\end{lem}

\begin{Proof} We can and will suppose $W$ is closed.
Since the (closed) balls form a basis for the topology, any open covering
can be refined to a covering form by closed balls. So we can and
will suppose we have a covering $\{U_i\}_{i\in I}$ of $W$ by closed
balls, and we have to show the covering has a finite subcovering.
Suppose it has not. Take $\delta\in \Tog$ topologically nilpotent.
For any $n\ge 1$, take inductively $B_n\in S_{\delta^n}$ a ball such
that $B_{n+1}\subset B_n$ and $\{U_i\cap B_n\}_{i\in I}$ has no
finite subcovering of $B_n\cap W$. Since $\lim_{n\to \infty}
\varrho(B_n)=0$ and they are nested balls, they intersection is a
point $p$. Since $B_n\cap W\ne \emptyset$ and $W$ is closed, $p\in
W$. Let $i_0\in I$ be such that $p\in U_{i_0}$. Then
$\varrho(U_{i_0})\ge \delta^{n_0}$ for some $n_0$, hence
$U_{i_0}\supset B_{n_0}$ since both have center $p$, hence
$U_{i_0}\cap B_{n_0}=B_{n_0}$ alone is a finite subcovering of
$B_{n_0}\cap W$, which contradicts the construction of the $B_n$.
\end{Proof}

In order to show the theorem, we will introduce some notations. We
denote by $S:=\{\vp_1,\dots,\vp_g,\vp_1^{-1},\dots,\vp_g^{-1}\}$.
Given $\psi\in S$, we denote
$$B(\psi):=\begin{cases}
B_i \mbox{ if } \psi=\vp_i^{-1}\\
B_{i+g} \mbox{ if } \psi=\vp_i
\end{cases}$$
and $F:=\bigcap_{\psi\in S} \psi^{-1}(B(\psi))$, which we call the
fundamental domain associated to $S$ and the balls $B_i$. Note that
$\infty \in \psi^{-1}(B(\psi))\supset B(\psi)$, so $\infty\in F$,
hence it is not empty. Note also that
$\vpB{\psi}{B(\psi^{-1})}=B(\psi)$ by definition.

Given an element $w\in \Gamma$, $w\ne \id$, a reduced word
expression for $w$ is an expression as a product
$w=\psi_s\cdots\psi_1$ where $\psi_i\in S$ for all $i$ and
$\psi_{i-1}\ne \psi_{i}^{-1}$ for all $i\ge 1$.

\begin{lem}\label{free} \begin{enumerate}\item For any
$\psi,\phi\in S$, $\psi(B(\phi))\subsetneq B(\psi)$ if and only if
$\psi\ne \phi^{-1}$.
\item  Given $w=\psi_s\cdots\psi_1$ written as reduced
word, we have $w(F)\subset B(\psi_s)$.
\item As a consequence, $\Gamma$ is free with free generators
$\vp_1,\dots,\vp_g$.
\end{enumerate}
\end{lem}
\begin{Proof}
For any $\psi,\phi\in S$, $B(\psi)\cap B(\phi)=\emptyset$ if
$\phi\ne \psi$ and $\bP^1(K)\setminus \psi^{-1}(B(\psi)) \subset
B(\psi^{-1})$, so $B(\phi)\subset \psi^{-1}(B(\psi))$ if $\phi\ne
\psi^{-1}$, which shows the first part.

 We show the second part by induction on $s$. The case $s=1$
is clear since $F\subset\psi^{-1}B(\psi)$ for all $\psi\in S$, while
the induction step is done using the first part.

The final assertion now is clear: if the group is not free, then we
could write $\id$ as a non-trivial reduced word. But we have just
showed that any non-empty reduced word acts non-trivially.
\end{Proof}

Given the lemma, we know that the expression as reduced word of an
element of $\Gamma$ is unique. We will denote the length of such a
word (which depend on the generators we have) as $\ell(w)$. If
$w=\psi_1\cdots\psi_s=w'\psi_s$ is a reduced word, we denote the
ball $B(w):=\vpB{w'}{B(\psi_s)}$. Denote also by
$$\rho:=\min\{\rho_{i,j} \ | \ i\ne j\in\{1,\dots,2g\}\}.$$ We have
then that $d(B(\psi),B(\phi))\ge \rho^{-1}$ for all $\psi\ne\phi\in
S$.

\begin{lem}\label{B(w)} For any $w \in \Gamma$, $w\ne \id$, written
as reduced word as $w=\psi_1\cdots\psi_s$, consider
$w_i=\psi_1\cdots\psi_i$ for $i=1,\dots,s$, and $u_i=w_i^{-1}w$ for
$i=1,\dots,s-1$. Then for all $i<s$, $B(w)=w_iB(u_i)\subset B(w_i)$
and $d(B(w),B(w_i))^{-1}<\rho^{s-i}$.

Hence $\varrho(B(w))<\rho^{s-1}\varrho(B(\psi_1))$.
\end{lem}

\begin{Proof} The equality $B(w)=w_iB(u_i)$ is clear once we know
$w_2(B(\psi_s))$ is a ball. But this is easily deduced from the
previous lemma \ref{free} (1).

From the same lemma, we have $\psi_iB(\psi_{i+1})\subset B(\psi_i)$
for all $i<s$, hence
$$B(w_{i+1})=w_iB(\psi_{i+1})=w_{i-1}\psi_iB(\psi_{i+1})\subset
w_{i-1}B(\psi_i)=B(w_i).$$ Hence, by lemma \ref{pairofballs},
$$d(B(w_{i+1}),B(w_i))=d(B(\psi_{i+1}),\vpB{{\psi_i^{-1}}}{B(\psi_i)})=
d(B(\psi_{i+1}),B(\psi_i^{-1}))\ge \rho^{-1}.$$

Hence
$$d(B(w),B(w_i))=\prod_{j=i}^{s-1} d(B(w_{j+1}),B(w_j))\ge
\rho^{i-s}.$$ The final assertion is clear from this inequality,
since $d(B(w),B(w_1))=\varrho(B(w_1))\varrho(B(w))^{-1}$.
\end{Proof}

As a consequence, we get that for any $\delta\in \Tog$, the set of
$w\in\Gamma$ such that $\varrho(B(w))>\delta$ is finite.

\begin{cor}\label{fundamentaldomain} For any $w\in \Gamma$, $w\ne
\id$, $w(F)\cap F = \emptyset$ if and only if $w\notin S$, while
$$\psi(F)\cap F=\psi(B(\psi^{-1}))\cap B(\psi)$$ if
$\psi\in S$.
\end{cor}
\begin{Proof} First note that for any $\psi\in S$,
$$B(\psi)\cap F=\bigcap_{\vp\in S} B(\psi)\cap \vp^{-1}B(\vp))=B(\psi)\cap \psi
B(\psi^{-1})$$ since $$ B(\psi)\cap \vp^{-1}B(\vp))=
\begin{cases} B(\psi) & \text{if } \vp\ne \psi^{-1}\\
B(\psi)\cap \psi B(\psi^{-1}) & \text{if } \vp=
\psi^{-1}.\end{cases}$$

Now, note that $$\psi(F)=\bigcap_{\vp} \psi\vp^{-1}B(\vp) \subset
B(\psi)$$ for any $\psi\in S$, so $$\psi(F)\cap F\subset B(\psi)\cap
F=B(\psi)\cap \psi B(\psi^{-1}).$$ But the reverse inclusion is
easy.

Finally, suppose $w\in \Gamma\setminus S$, $w\ne \id$, so
$\ell(w)\ge 2$. Then $w=w'\psi$ for some $\psi\in S$, $w'\ne id$,
and $\ell(w)=\ell(w')+1$. Then $w(F)=w'(\psi(F))\subset
w'B(\psi)\subset B(w')$ by lemma \ref{B(w)}. But $w'B(\psi)\cap \vp
B(\vp^{-1})=\emptyset$ if $w'=\vp w''$ as a reduced word, since $w''
B(\psi)\cap B(\vp^{-1})=\emptyset$. This last equality is obvious if
$w''=id$, since $\psi\ne \vp^{-1}$. Finally, if $w''=\tau w'''$ as
reduced word, it is deduced from $w''B(\psi)\subset B(\tau)$ by
lemma \ref{B(w)} since $\tau\ne \vp^{-1}$.
\end{Proof}

\begin{pf1}
First we show that all elements $w\in \Gamma$, $w\ne \id$, are
hyperbolic. First of all, note that there exists $u,v\in\Gamma$,
$u\ne \id$, such that $w=vuv^{-1}$ and $\ell(u^2)=2\ell(u)$. Clearly
$w$ is hyperbolic if and only if $u$ it is by Lemma \ref{ppp}. Now,
$B(u^2)=uB(u)\subset B(u)$ by lemma \ref{B(w)}, and
$d(B(u),uB(u))^{-1}\le \rho^{\ell(u)}$, so it is topologically
hyperbolic. Now Proposition \ref{hyperbolicballs} implies that $u$
is hyperbolic.

To finish we only need to show that the orbit $\overline{\Gamma P}$
is compact for any $P\in \bP^1(K)$. We show this by proving that for
any point $P\in \bP^1(K)$ and for any $\delta$, $\Gamma P\backslash
\bigcup_{w\in \Gamma_{\delta}} B(w)$ is finite, where
$$\Gamma_{\delta}:=\{w'\psi\in \Gamma | \varrho(B(w'))>\delta,\ \varrho(B(w'\psi))\le \delta \mbox{ and } \ell(w'\psi)=\ell(w')+1\}$$
if $\delta<\max\{\varrho(B(w)) \ | \ \ell(w)=2\}$, and
$\Gamma_{\delta}=S$ if not. It is clear that $\Gamma_{\delta}$ is a
finite set, since we have a finite-to-one map to a subset of the
$w'\in\Gamma$ such that $\varrho(B(w'))>\delta$, which is finite by
lemma \ref{B(w)}. The proof is then finished by using Lemma
\ref{Heine-Borel}.

To show the assertion, we can and will suppose that
$\delta<\max\{\varrho(B(w)) \ | \ \ell(w)=2\}$. Observe that then
for any $\tau\in \Gamma_{\delta}$, $B(\tau)\subset \bigcup_{w\in
\Gamma_{\delta}} B(w)$.

If $P\in F$, for any $w\in \Gamma$, $w\ne \id$, then $w(P)\in B(w)$.
Hence $$\Gamma P\backslash \bigcup_{w\in \Gamma_{\delta}}
B(w)\subset \{\tau(P)\ | \ \varrho(B(\tau))>\delta\}$$ which is
clearly finite. The same is true if $P\in \vp(F)$ for some $\vp\in
\Gamma$, as they have the same $\Gamma$-orbits.

Finally, if $P\in \bP^1(K)\setminus \bigcup_{\vp\in \Gamma} \vp(F)$,
then $P\in B(w)$ for some $w\in \Gamma_{\delta}$. Therefore
$\vp(P)\in \bigcup_{w\in \Gamma_{\delta}} B(w)$ for all $\vp\in
\Gamma$.
\end{pf1}

In order to show a reverse theorem, asserting that any Schottky
group $\Gamma$ with $\infty \notin \cL_{\Gamma}$ has potentially a
good fundamental domain, we need an elementary result on graphs with
$\Tog$-weights. For a given totally ordered (multiplicative) group
$\Tog$, a graph with $\Tog$-weights is a graph $G$ together with a
map $w:E(G)\to \Tog_{>1}$ assigning to any edge $e$ an element
$w(e)\in \Tog_{>1}$. For a given $c=\sum_{e\in I} e$, where $I$ is a
finite set, one defines $w(c)=\prod_{e\in I}w(e)$.

\begin{lem}\label{goodedges} Suppose $(G,w)$ is a finite graph with $\Tog$-weights such
that, for any cycle $c\in H_1(G,\bZ)$, $w(c)^{-1}$ is topologically
nilpotent. Let $g=\rank_{\bZ}(G)$ be the genus of $G$. Then there
exists edges $e_1,\dots,e_g$ such that $w(e_i)^{-1}$ is
topologically nilpotent for all $i=1,\dots,g$ and
$G\setminus\{e_1,\dots,e_g\}$ is a (spanning) tree.
\end{lem}

\begin{Proof}  By induction on $g$. If $g=0$ there is nothing to prove.
Let $c_1=\sum_{e\in I_1}e$ be a cycle. Then there exists $ e_1\in
I_1$ such that $w(e_1)^{-1}$ is topologically nilpotent, since
otherwise $w(c_1)^{-1}=\prod_{e\in I_1}w(e)^{-1}$ would be a product
on non-topologically nilpotent elements, so non-topologically
nilpotent. But then the graph $G_1:=G\setminus\{e_1\}$ has genus
$g-1$ and verifies the hypothesis of the lemma. \end{Proof}

\begin{teo}\label{reversetheorem} Let $\Gamma$ be a Schottky group in $\PGL(K)$ of rank $g$
such that $\infty \notin \cL_{\Gamma}$. Then there exists an
extension $L/K$ of type $(2,\dots,2)$, closed balls
$B_1,\dots,B_{2g}$  in $\bP^1(L)$ with $d(B_i,B_j)^{-1}$
topologically nilpotent for all $i\ne j\in\{1,\dots,g\}$ and
generators $\vp_1,\dots,\vp_g$ of $\Gamma$ such that
$\vpB{\vp_i}{B_i}=B_{i+g}$ and $\infty\in\vp_i(B_i)$ for
$1\in\{1,\dots,g\}$.
\end{teo}

\begin{Proof} Consider the finite graph $G_{\Gamma}$ constructed in
section 7. We assign weights to any edge $e$ of $G_{\Gamma}$ by
lifting them to $\cT_{\Gamma}$ and assigning them the distance
between the two extreme vertices. Any cycle $c$ in $G_{\Gamma}$
correspond to the projection of a path between a vertex $v$ and
$\vp(v)$, for some $\vp\in \Gamma$. Then
$d(v,\vp(v))^{-1}=w(c)^{-1}$ is then topologically nilpotent by
proposition \ref{hyperbolicballs}. Hence we are under the hypothesis
of the lemma \ref{goodedges}, so we can take edges $e_1,\dots,e_g$
such that $T:=G_{\Gamma}\setminus\{e_1,\dots,e_g\}$ is a spanning
tree and with $w(e_i)^{-1}$ topologically nilpotent for any
$i=1,\dots,g$.

Since $\cL_{\Gamma}$ is compact and $\infty\notin \cL_{\Gamma}$,
$\cL_{\Gamma}$ is contained in the union of the closed balls
corresponding to vertices of $\cT_{\Gamma}$, hence in the union of a
finite number of them. Two options can happen: either there is one
of such balls that contains all $\cL_{\Gamma}$, or there are two of
them joined by an edge. This is because, if there are three disjoint
balls corresponding to vertices of $\cT_{\Gamma}$, there is a vertex
of $\cT_{\Gamma}$ whose corresponding ball contain two of them. We
call this edge or this vertex in either case the root of
$\cT_{\Gamma}$.

We lift the tree $T$ to a subtree $\Upsilon\subset \cT_\Gamma$
containing the root of $\cT_{\Gamma}$.

The edges $e_i$ can be lifted in two distinct ways to edges in
$\cT_{\Gamma}$ attached to $\Upsilon$; denote them by $\bar e_i$ and
$\bar e_{i+g}$ respectively. There are then elements $\vp_i\in
\Gamma$ sending $\vp_i(\bar e_i)=\bar e_{i+g}$ (and the vertex $v_i$
touching $\Upsilon$ to $v_{i+g}$), which generate the group $\Gamma$
(being $\Gamma$ the fundamental group of $G_{\Gamma}$ and
$T_{\Gamma}$ its universal covering).

Now, given an edge $\bar e$, seen as a subset of the $\Tog$-tree
$\BK$, we consider the middle point $B_e$ in $\cT_{K'}$, where
$K'/K$ is an extension of degree at most 2; if $\bar
e=[B(p,\delta),B(p,\rho)]$, then $B_e=B(p,\sqrt{\delta\rho})$, and
this is defined in the totally ordered group $\Tog'$ associated to
the splitting field $K'$ of the polynomial $X^2-\alpha$ for some
$\alpha\in K$ with $|\alpha|=\delta\rho$. And all the edges $\bar
e_i$ are of this form because their extremes are contained in the
root or in one of the extremes of the root. The corresponding balls
$B_i:=B_{e_i}$ are then pairwise disjoint, and
$\vp{{\vp_i}}{B_i}=B_{i+g}$. Clearly $\infty \in \vp_i(B_i)$.

Now, the distance $d(B_i,B_j)$ for $i\ne j$ is equal to length of
the path from $B_i$ to $B_j$. This path can be decomposed as the
(half)-edge from $B_i$ to $v_i$, the path from $v_i$ to $v_j$ and
the edge from $v_j$ to $B_i$. But the length of the edge from $B_i$
to $v_i$ is equal to $\sqrt{w(e_i)}=\sqrt{\varrho(\vp_i)}$, with
inverse topologically nilpotent. So a fortiori the same is true for
$d(B_i,B_j)$.
\end{Proof}

A similar argument can be used to show the following result.

\begin{cor}\label{GraphtoGroup} Let $K$ be a field complete with respect to a
$\Tog$-valuation. Let $G$ be a finite graph with $\Tog$-weights and
with no vertices of valence $\le 2$. Suppose that $w(c)^{-1}$ is
topologically nilpotent for any cycle $c$. Then there exists a
finite extension $L/K$ and a Schottky group $\Gamma\subset
\PGL_2(L)$ with $G_{\Gamma}\cong G$ as graph with $\Tog$-weights.

Moreover, if the residue field $k$ of $K$ is not finite or the
valence of all vertices of $G$ is not larger than $\#k+1$ with at
least one vertex of valence $\le \#k$, we can take $L=K$.
\end{cor}

To prove the corollary we will used the following elementary lemma.

\begin{lem}\label{starinT} Let $B$ be a closed ball in $K$ and take
$\rho_1,\dots,\rho_n\in \Tog_{>1}$. Suppose that the residue field
$k$ of the valuation has more than $n$ elements. Then there exists
balls $B_1,\dots,B_n \subset B$, pairwise disjoint, and with
$d(B,B_i)=\rho_i$ for all $i\in\{1,\dots,n\}$.
\end{lem}

\begin{Proof} If $p_1$ and $p_2$ are points in $K$ such that $B=t(p_1,p_2,\infty)$,
consider the polynomial automorphism $\varphi$ sending $p_1$ to 0
and $p_2$ to $1$ (and hence $\infty$ is fixed). Then the ball $B$ is
send to the ball $\cO=B(0,1)$, and we are reduced to this case.
Consider $\delta:=\min\{\rho_i^{-1}\ | \ i\in\{0,\dots,n\}\}$ and
$B(0,\delta):=\{x\in K \ | \ |x|\le \delta \}$, which is an ideal of
$\cO$. Then the natural map $\psi_{\rho}\colon \cO\to
\cO/B(0,\delta)$ send points in $\cO$ to the same image if and only
if they are contained in a common ball of radius $\le \delta$. Hence
we can find $n$ disjoint balls in $\cO$ of radius $\rho_i^{-1}$ for
$i\in\{1,\dots,n\}$, as long as $n$ is less than or equal to the
number of elements in $\cO/B(0,\delta)$, which is larger or equal to
the number of elements in $k=\cO/\mm$.
\end{Proof}

\begin{pf2} Given such a graph $G$ of genus $g\ge 2$, we use lemma \ref{goodedges} to
find $\{e_1,\dots,e_g\}$ such that $T:=G\setminus\{e_1,\dots,e_g\}$
is a spanning tree and $w(e_i)^{-1}$ is topologically nilpotent for
all $i\in \{1,\dots,g\}$. We choose a ball $B$ in $K$, which will be
the root of $\Gamma$, and a vertex $v$ in $T$ with valence $\le
\#k$. Using the previous lemma \ref{starinT}, we lift the start in
$T$ of the vertex $v$ by choosing balls inside $B$, pairwise
disjoint, and with distance the length of the corresponding edges of
$T$ with one extreme in $v$ (defined in a suitable extension $L$ if
the field $k$ is finite and too small). For any one of this new
balls, we repeat the process the any of the vertices corresponding
to $T\setminus\{v\}$. At the end we lifted the tree $T$ to a subtree
$\Upsilon\subset \cT_\Gamma$ with root in $B$. We do the same
process like in the proof of theorem \ref{reversetheorem} to attach
some half edges to $\Upsilon$ corresponding to liftings of the $e_i$
to get a configuration of balls verifying the hypothesis of theorem
\ref{fundamental}, and hence a Schottky group $\Gamma$, defined in
an extension $L'/L$ of type $(2,\dots,2)$. By construction
$G_{\Gamma}\cong G$. It only remains to show that the group $\Gamma$
can be chosen already defined over $L$. This is done by proving that
the hyperbolic automorphisms constructed in the lemma \ref{twoballs}
in our case can be chosen in $K$, since the distance between the
balls is in $\Tog$, and not in the totally ordered group
corresponding to $L'$.
\end{pf2}

We believe that the condition asking for at least one vertex with
valence $\le \#k$ is not necessary, but we would need a notion of
good fundamental allowing $\infty\in \cL_{\Gamma}$.

One can use this result to give a new proof of a known result of B.
Conrad in appendix B in \cite{baker}. The idea behind this proof is
probably well-known, as it can be deduced of the already known
result.

\begin{cor} Let $R$ be a complete discrete valuation domain, with
field of fractions $K$ and residue field $k$. Let $G$ be a graph
with no vertices of valence $1$, and with the valences of all
vertices $\le \#k+1$, and one vertex of valence $\le \#k$. Then
there exists a smooth projective curve over $K$ with a split
semistable reduction and dual graph of the reduction isomorphic to
$G$.
\end{cor}

\begin{Proof} If the genus of $G$ is 1, then the graph is a cycle
graph, and the Tate elliptic curve $E=\mathbb{G}m/q^{\bZ}$ with
valuation equal to the number of vertices solves the problem. If the
genus is $\le 2$, we take the $\bZ$-graph associated by deleting the
vertices of valence 2 and with the weight of the new edges the
number of old edges that contain. Applying the corollary
\ref{GraphtoGroup} we get a Schottky group $\Gamma$ with associated
graph isomorphic to $G$. Associated to $\Gamma$ one constructs the
Mumford curve obtained as analytic curve as the quotient
$(\bP^1\setminus \cL_{\Gamma})/\Gamma$. It has a split semiestable
model with graph of the reduction isomorphic to $G_{\Gamma}$ (see
for example \cite{mumford1972analytic}, page 163, or \cite{GvdP},
remark III.2.12.3).
\end{Proof}

Note that Conrad's result is stronger in the sense that he shows the
existence of a smooth and projective curve with reduction graph any
given graph, not just graphs without valence one vertices. However,
he proofs the result only when $k$ is an infinite field (but his
proof can be probably adapted to show the result above).

\bibliographystyle{plain}

\begin{thebibliography}{10}

\bibitem{baker}
Matt Baker.
\newblock Specialization of linear systems from curves to graphs.
\newblock {\em Algebra \& Number Theory}, 2(6):613--653, 2008.

\bibitem{berkovich2012spectral}
Vladimir~G Berkovich.
\newblock {\em Spectral theory and analytic geometry over non-Archimedean
  fields}.
\newblock Number~33. American Mathematical Soc., 2012.

\bibitem{bruhat1972groupes}
Fran{\c{c}}ois Bruhat and Jacques Tits.
\newblock Groupes r{\'e}ductifs sur un corps local.
\newblock {\em Publications Math{\'e}matiques de l'IH{\'E}S}, 41(1):5--251,
  1972.

\bibitem{brum}
G.~W. Brumfiel.
\newblock The tree of a non-{A}rchimedean hyperbolic plane.
\newblock In {\em Geometry of group representations ({B}oulder, {CO}, 1987)},
  volume~74 of {\em Contemp. Math.}, pages 83--106. Amer. Math. Soc.,
  Providence, RI, 1988.

\bibitem{chiswell}
I.~Chiswell.
\newblock {\em Introduction to $\Lambda$-trees}.
\newblock World Scientific, 2001.

\bibitem{EnglerPrestel}
Antonio~J. Engler and Alexander Prestel.
\newblock {\em Valued fields}.
\newblock Springer Monographs in Mathematics. Springer-Verlag, Berlin, 2005.

\bibitem{GvdP}
Lothar Gerritzen and Marius van~der Put.
\newblock {\em Schottky groups and Mumford curves}.
\newblock Springer-Verlag, 1980.

\bibitem{HuEt}
Roland Huber.
\newblock {\em \'Etale cohomology of rigid analytic varieties and adic spaces}.
\newblock Aspects of Mathematics, E30. Friedr. Vieweg \& Sohn, Braunschweig,
  1996.

\bibitem{MoSh}
John~W. Morgan and Peter~B. Shalen.
\newblock Valuations, trees, and degenerations of hyperbolic structures. {I}.
\newblock {\em Ann. of Math. (2)}, 120(3):401--476, 1984.

\bibitem{mumford1972analytic}
David Mumford.
\newblock An analytic construction of degenerating curves over complete local
  rings.
\newblock {\em Compositio Mathematica}, 24(2):129--174, 1972.

\bibitem{Serretrees}
Jean-Pierre Serre.
\newblock {\em Trees}.
\newblock Springer-Verlag, 1980.

\bibitem{GVX}
Iago~Gin\'e V\'azquez and Xavier Xarles.
\newblock The abel-jacobi map for mumford curves via integration.
\newblock {\em Preprint}, 2016.

\end{thebibliography}

\end{document}